\documentclass[12pt]{article}
\usepackage{amsmath,amsthm,amssymb}
\usepackage{mathtools}

\renewcommand{\theequation}{\mbox{\arabic{section}.\arabic{equation}}}

\newtheorem{theorem}{Theorem}[section] 
\newtheorem{lemma}{Lemma}[section]
\newtheorem{cor}{Corollary}[section]
\newtheorem{remark}{Remark}[section]

\newtheorem{proposition}{Proposition}[section]

\numberwithin{equation}{section}

\newcommand{\D}{{\rm d}}
\newcommand{\dx}{\, \D x}
\newcommand{\dy}{\, \D y}
\newcommand{\dz}{\, \D z}
\newcommand{\ds}{\, \D s}
\newcommand{\dr}{\, \D r}

\newcommand{\dis}{\displaystyle}

\newcommand{\msp}{\;\;}
\newcommand{\fsp}{\quad\;}
\newcommand{\psp}{\,}

\newcommand{\rz}{\mathbb{R}}
\newcommand{\nz}{\mathbb{N}}

\newcommand{\eps}{\varepsilon}
\newcommand{\iom}{\int_{\Omega}}

\newcommand{\klauf}{\left(\begin{array}}
\newcommand{\klzu}{\end{array}\right)}

\newcommand{\udel}{u_{\delta}}
\newcommand{\gam}[2]{\Gamma_{#1,\delta}^{#2}}

\newcommand{\hatw}{\hat{w}}

\newcommand{\aep}{A_{i,\bar{\eps}}}
\newcommand{\ovu}{\overline{u}}

\newcommand{\csob}{\mathcal{C}_{{\rm Sob}}}
\newcommand{\cbv}{\mathcal{C}_{{\rm BV}}}

\newcommand{\wnulla}{{\stackrel{\circ}{W}}\!\prescript{1}{A}(\Omega)}
\newcommand{\wnull}{{\stackrel{\circ}{W}}\!\prescript{1}{1}(\Omega)}

\newcommand{\wnullaa}{{\stackrel{\circ}{W}}\!\prescript{1\mbox{\,\, }}{1,A}(\Omega)}

\title{Variational problems of splitting-type with mixed linear- superlinear growth conditions}
\author{Michael Bildhauer \& Martin Fuchs}
\date{}

\newcommand{\reff}[1]{(\ref{#1})}

\usepackage{hyperref}

\hypersetup{
   pdfnewwindow=true,
   colorlinks=false,
   linkcolor=black,
   citecolor=green,
   filecolor=magenta,
   urlcolor=cyan,
}

\begin{document}

\parindent0em
\maketitle

\newcommand{\op}[1]{\operatorname{#1}}
\newcommand{\bv}{\op{BV}}
\newcommand{\mub}{\overline{\mu}}
\newcommand{\muhat}{\hat{\mu}}

\newcommand{\hypref}[2]{\hyperref[#2]{#1 \ref*{#2}}}
\newcommand{\hypreff}[1]{\hyperref[#1]{(\ref*{#1})}}

\newcommand{\ob}[1]{^{(#1)}}

\begin{abstract}
Variational problems of splitting-type with mixed linear- superlinear growth conditions are considered. In the twodimensional case the 
minimizing problem is given by
\[
J [w] = \iom \Big[f_1\big(\partial_1 w\big) + f_2\big(\partial_2 w\big)\Big] \dx \to \min
\]
w.r.t.~a suitable class of comparison functions. Here $f_1$ is supposed to be a convex energy density with linear growth,
$f_2$ is supposed to be of superlinear growth, for instance to be given by a $N$-function or just bounded from below by a $N$-function.
One motivation for this kind of problem located between the well known splitting-type problems of superlinear growth and the
splitting-type problems with linear growth (recently considered in \cite{BF:2020_3}) is the link to mathematical problems in plasticity 
(compare \cite{Te:2018_1}). Here we prove results on the appropriate way of relaxation including approximation procedures, duality,
existence and uniqueness of solutions as well as some new
higher integrability results.\footnote{AMS-Classification: 49J45, 49N60}
\end{abstract}

\parindent0ex
\section{Introduction}\label{intro}

In the last decades the study of variational problems with nonstandard growth conditions 
developed to one of the main topics in the calculus of variations and related areas. 
We do not want to go into details and will not present the historical line of this development. 
The reader will find this background information, for instance,
in the recent paper \cite{BM:2020_1}.\\

Let us just mention a few aspects, which serve as a motivation for the manuscript at hand.\\

As one of the first main contributions Giaquinta considered the most common prototype of 
energies with $(p,q)$-growth in the sense of minimizing a 
splitting functional
\begin{equation}\label{hist 1}
\iom f(\nabla u) \dx = \iom \Big[f_1\big(\partial_1 u\big) +f_2(\partial_2 u\big)\Big] \dx \to \min 
\end{equation}
in a suitable class of comparison functions, where $f_1$ and $f_2$ are supposed to have different growth rates larger than $1$. 
In \cite{Gi:1987_1} he presented a famous counterexample which shows
that in general and even in the scalar case we cannot expect the smoothness of solutions to this variational problem.\\

Of course \reff{hist 1} also serves as a motivation to study variational problems with non-uniform ellipticity conditions.
As one variant we may consider energy densities of class $C^2$ satisfying with different exponents $1 < p < q$
($c_1$, $c_2 >0$, $\xi$, $\eta \in \rz^{nN}$)
\begin{equation}\label{hist 2}
c_1 \big(1+|\xi|^2\big)^{\frac{p-2}{2}} |\eta|^2 \leq D^2f(\xi)(\eta,\eta) \leq c_2 \big(1+|\xi|^2\big)^{\frac{q-2}{2}} |\eta|^2 \psp .
\end{equation}

Here a lot of important contributions in the scalar and also in the vectorial setting can be found, we just mention \cite{Ma:1991_1} as one
central reference in the long series of papers in this direction.\\

Related to \reff{hist 2}, Frehse and Seregin (\cite{FrS:1999_1}) considered plastic materials with logarithmic hardening, i.e.~the energy density
$f(\xi) =|\xi| \ln\big(1+|\xi|\big)$ of nearly linear growth. Due to \cite{MS:1999_1} we have full regularity for this particular kind of model. \\

We finally pass to the case of linear growth problems for which the energy density is of (uniform) linear growth
w.r.t.~the gradient but just satisfies a non-uniform ellipticity condition in the sense of \reff{hist 2}. 
Of course the minimal surface case is the most prominent representative for this kind of problems
satisfying with suitable constants $a_1$, $b_1$, $c_1$, $c_2 >0$, $a_2$, $b_2 \geq 0$, for all $\xi$, $\eta \in \rz^{nN}$
and with some exponent $\mu >1$
\begin{equation}\label{hist 3}
\begin{array}{rcccl}
a_1 |\xi| - a_2 &\leq&f\big(\xi\big) &\leq & b_1 |\xi| + b_2\psp ,\\[2ex]
c_1 \big(1+|\xi|^2\big)^{-\frac{\mu}{2}}|\eta|^2  &\leq & D^2f(\xi)(\eta,\eta)& \leq &
c_2 \big(1+|\xi|^2\big)^{-\frac{1}{2}} |\eta|^2 \psp .
\end{array}
\end{equation}

In the minimal surface case we have $\mu =3$ and we like to mention the pioneering work of Giaquinta, Modica and Sou\v{c}ek
\cite{GMS:1979_1}, \cite{GMS:1979_2} in the list of outstanding contributions. In \cite{GMS:1979_1} and \cite{GMS:1979_2} 
a suitable relaxation is discussed together with a subsequent proof of apriori estimates. We note that the uniqueness of solutions
in general is lost by passing to the relaxed problem.\\

In \cite{BF:2003_1}, condition \reff{hist 3} was introduced defining a class of $\mu$-elliptic energy densities. The regularity theory for minimizers
was studied in a series of subsequent papers, compare, e.g., \cite{Bi:2002_1}. We also like to mention the Lipschitz estimates of 
Marcellini and Papi \cite{MP:2006_1}, which cover a broad class of the functionals we discussed up to now.\\

Very recently, the authors \cite{BF:2020_3} considered variational problems of splitting-type as given in \reff{hist 1} but now
with two energy parts being of linear growth. Here it turns out that the right-hand side of the ellipitcity condition in \reff{hist 3}
is no longer valid and we just have for all $\xi$, $\eta \in \rz^{nN}$ with a positive \mbox{constant $c$}
\[
D^2f(\xi)(\eta,\eta) \leq c |\eta|^2 \psp .
\] 
Nevertheless, some natural assumptions still imply
regularity and uniqueness properties of solutions to the relaxed problem.\\

In the manuscript at hand we follow this line of studying variational problems of splitting-type by now considering
variational problems with mixed linear- superlinear growth conditions. \\

A first step in this direction was already made in Chapter 6 of \cite{Bi:1818}. The results given there follow from suitable apriori estimates which
are available if the non-uniform ellipticty is not too bad. 
This leads to the analysis of the set of cluster points of minimizing sequences
and to some kind of local interpretation for the stress tensor, 
although the existence and the uniqueness of dual solutions were not established (compare Remark 6.15 of \cite{Bi:1818}).\\

We like to finish this introductory remarks by mentioning a prominent application of a mixed linear- superlinear growth problem, which 
in \cite{Te:2018_1} is discussed as the Hencky plasticity model. This problem takes the form (compare (4.17), Chapter I, of \cite{Te:2018_1})
\[
\inf_{v\in \mathcal{C}_a}\Bigg\{ \iom (\op{div} v)^2\dx + \iom \psi\big(\eps^D(v)\big) \dx - L(v)\Big\}
\]
with a suitable class $\mathcal{C}_a \subset W^{1,2}(\Omega;\rz^3)$, some volume force $L$ and the deviatoric part $\eps^D(v)$ of the symmetric
gradient $\eps(v)$. Here $\op{div} v$ enters with quadratic growth while the function $\psi$ is of linear growth
w.r.t.~the tensor $\eps^D(v)$.
Let us note that this plasticity model is based on the dual point of view, i.e.~the so-called sur-potential $\psi$ is introduced
through the conjugate function $\psi^*$. In general a more explicit expression cannot be given (see Remark 4.1, p.~75 of \cite{Te:2018_1}).
One particular example is of the form $\psi(\xi^D) = \Phi\big(|\xi^D|\big)$ with (for some positive constants $k$, $\nu$)
\[
\Phi(s) = \left\{\begin{array}{lcr}
\nu s^2&\mbox{if}&|s| \leq \frac{k}{\sqrt{2}\nu} \psp ,\\[2ex]
\sqrt{2} k |s| - \frac{k^2}{2\nu}&\mbox{if}& |s| \geq \frac{k}{\sqrt{2} \nu}\psp .
\end{array}\right.
\]

However, at this stage our main difficulty in comparison to the known results for the Hencky model is quite hidden. 
We postpone a refined discussion to Remark \ref{remark newapprox 1}.\\
 
Now let us introduce the general framework of our considerations in a more precise way. For the sake of notational simplicity we
restrict our considerations to the case that a linear growth condition is satisfied in only one coordinate direction.\\

Suppose that $\Omega \subset \rz^n$ is a bounded Lipschitz domain and let $f$: $\rz^n \to \rz$
be an energy density of class $C^2(\rz^n)$ which is decomposed in the form
\begin{equation}\label{intro not 1}
f(\xi) =f_1 (\xi_1) + f_2(\xi_2)\psp , \fsp \xi =(\xi_1,\xi_2) \in \rz \times \rz^{n-1} \psp .
\end{equation}

Here we assume that $f_1$: $\rz \to \rz$ and $f_2$: $\rz^{n-1} \to \rz$ are convex functions of class 
$C^2(\rz)$ and $C^2(\rz^{n-1})$, respectively, satisfying with $a_i$, $b_i\geq 0$, $a_1$, $a_3$, $b_1 >0$ and with a $N$-function $A$: $\rz \to \rz$:
\begin{eqnarray}\label{intro not 2}
&\dis a_1 |\xi_1| - a_2\leq  f_1(\xi_1)\leq  a_3 |\xi_1| + a_4\psp , \fsp \xi_1 \in \rz\psp ,\nonumber\\[2ex]
&\dis  b_1 A(|\xi_2|) - b_2 \leq f_2(\xi_2) \psp , \fsp \xi_2\in \rz^{n-1}\psp .&
\end{eqnarray}

We wish to note that we work with energy densities of class $C^2$ just for notational simplicity. This hypothesis just enters Section \ref{apriori}.
It is easy to check that for example no differentiability assumptions are needed in Section \ref{gen}, whereas the
results of Section \ref{dual} hold for densities of class $C^1$.\\ 

Having the decomposition \reff{intro not 1} in mind, we now suppose $n=2$ throughout the rest of this manuscript.
This essentially clarifies the notation while the main ideas and results remain unchanged.\\

For the definition and the properties of $N$-functions and Orlicz-Sobolev space we refer to the monographs
\cite{RR:1991_1} or \cite{KR:1961_1}. The basics needed here are summarized in \cite{FO:1998_1} or \cite{FL:1998_1}.
We suppose that $A$: $[0,\infty) \to [0,\infty)$ satisfies
\[
\begin{array}{ll}
(N1)&\mbox{$A$ is continuous, strictly increasing and convex.}\\[2ex]
(N2)&\dis \lim_{t\downarrow 0} \frac{A(t)}{t}= 0 \fsp\mbox{and}\fsp \lim_{t\to \infty} \frac{A(t)}{t} =  \infty\psp .\\[4ex]
(N3)&\dis\mbox{There exist constants $k$, $t_{0} >0$ such that $A(2t)  \leq k A(t)$}\\ 
&\mbox{for all $t \geq t_{0}$}\psp , 
\end{array}
\]
where $(N3)$ is called a $\Delta_2$-condition near infinity. We note that there exists an exponent $p \geq 1$ such that with some constant
$c>0$
\begin{equation}\label{A growth}
c |t|^p \leq A(t) \fsp\mbox{for all}\msp t\gg 1\psp .
\end{equation}

Moreover we suppose some kind of triangle inequality for $f_2$: there exists a real number $c_3>0$ such that
for all $t$, $\hat{t}\in \rz$
\begin{equation}\label{intro not 5}
f_2 \big(t + \hat{t}\big) \leq c_3\Big[ f_2\big(t\big) + f_2\big(\hat{t}\big)\Big] \psp .
\end{equation}
This condition, for instance, follows from the convexity of $f_2$ together with some $\Delta_2$-condition.\\

For the definition of the Sobolev spaces $W^{k}_{p}$ and their local variants we refer to the textbook of Adams (\cite{Ad:1975_1}), 
the notation needed later in the case of functions of bounded variation can be found, e.g., in the monographs \cite{AFP:2000_1} and \cite{Gi:1984_1}.
For the sake of completeness we recall the definition of the Orlicz-Sobolev space generated by a $N$-function $A$ satisfying $(N1)$-$(N3)$
(see \cite{RR:1991_1}, \cite{KR:1961_1}).\\

In the following we suppose that the bounded Lipschitz domain $\Omega$ is normal w.r.t.~the $x_2$-axis (compare the
approximation arguments presented in \mbox{Section \ref{approx}}), i.e.~there exist Lipschitz functions $\kappa_1$, $\kappa_2$:
$(a,b)\to \rz$ such that 
\begin{equation}\label{intro not normal}
\Omega = \big\{x\in \rz^2:\psp x_1\in (a,b), \psp \kappa_1(x_1) < x_2 < \kappa_2(x_1) \big\} \psp .
\end{equation}
The space
\begin{eqnarray*}
L_A(\Omega) &:=& \Bigg\{u:\psp \Omega \to  \rz: \psp \mbox{$u$ is a measurable function such that}\\ 
&&\mbox{there exists $\lambda >0$ with} \iom A \big(\lambda |u|\big) \dx < +\infty\Bigg\} 
\end{eqnarray*}
is called Orlicz space equipped with the Luxemburg norm
\[
\|u\|_{L_A(\Omega)} = \inf\Bigg\{ l >0:\psp \iom A\Bigg(\frac{|u|}{l}\Bigg) \dx \leq 1\Bigg\}\psp .
\]
The Orlicz-Sobolev space is given by
\[
W^1_A(\Omega) = \Bigg\{u:\psp \Omega \to \rz: \psp \mbox{$u$ is a measurable function}, \psp 
u, \psp \big|\nabla u\big| \in L_A(\Omega)\Bigg\} 
\]
with norm
\[
\|u\|_{W^1_A(\Omega)} = \|u\|_{L_A(\Omega)} + \|\nabla u\|_{L_A(\Omega)}\psp .
\]

The closure in $W^1_A(\Omega)$ of $C^\infty_0(\Omega)$-functions w.r.t.~this norm is according to
Theorem 2.1 of \cite{FO:1998_1} (recall that we suppose $(N3)$)
\begin{equation}\label{intro not 6}
\wnulla =W^1_A(\Omega) \cap \wnull \psp .0
\end{equation}

We additionally use the notation
\begin{eqnarray*}
W^{1}_{1,A}(\Omega) &:=& \Big\{ w\in W^{1,1}(\Omega): \psp \partial_2 w \in L_A(\Omega)\Big\}\psp ,\\[2ex]
E[v] &:=& \iom f_2(v) \dx \psp , \fsp v \in L_A(\Omega) \psp .
\end{eqnarray*}

In formal accordance with  \reff{intro not 6} we define the class
\[
\wnullaa :=W^1_{1,A}(\Omega) \cap \wnull \psp .
\]
observing that this set is just the completion of $C^\infty_0(\Omega)$ in $W^1_{1,A}(\Omega)$ w.r.t.~the natural
norm of $W^1_{1,A}(\Omega)$.\\

The main classes of functions under consideration are:
\begin{eqnarray*}
\csob &:=& u_0+\wnullaa \psp ,\\[2ex]
\cbv &:=&\Bigg\{ w \in \bv(\Omega):\psp \| \partial_2 w\|_{L_A(\Omega)} < \infty ,
\psp  (w-u_0)\nu_2=0 \msp \mbox{$\mathcal{H}^1$-a.e.~on $\partial \Omega$}\Bigg\} \psp ,
\end{eqnarray*}
where the boundary values $u_0$ are always supposed to be of class $W^{1}_{\infty}(\Omega)$ (see also Remark 6.12 in \cite{Bi:1818}).
Note that in the definition of the space $\cbv$ we require that the distributional derivative $\partial_2 w$ is generated by a function
from the space $L_{A}(\Omega)$. Moreover, we consider the $\bv$-trace of $w$ and denote by 
$\nu =(\nu_1,\nu_2)$ the outward unit normal to $\partial \Omega$. \\

With respect to these classes we consider the minimization problem
\begin{equation}\label{probsob}
J[w]:= \iom f(\nabla w) \dx \to \min\fsp\mbox{in the class}\msp \csob
\end{equation}

and its relaxed version
\begin{eqnarray}\label{probbv}
K[w]&:=& \iom f_1 \big(\partial_1^a w\big)\dx + \iom f_1^{\infty}\Bigg( \frac{\partial_1^s w}{|\partial_1^s w|}\Bigg)
\D |\partial_1^s w|\nonumber\\[2ex]
&& + \int_{\partial \Omega} f_1^\infty \big((u_0-w)\nu_1\big)\D \mathcal{H}^1
+ E\big[\partial_2 w\big] \nonumber\\[2ex] 
&=:&K_1[w] + E\big[\partial_2 w\big]  \to \min \fsp \mbox{in the class}\msp \mathcal{\cbv}\psp . 
\end{eqnarray}
Here $\nabla^a w$ denotes the absolutely continuous part of $\nabla w$ w.r.t.~the Lebesgue measure, 
$\nabla^sw$ represents the singular part.\\

As a matter of fact, problem \reff{probsob} in general is not solvable and one has to pass to the relaxed
version in order to have the existence of at least generalized minimizers. The approach via relaxation in the case of linear growth
is well-known and outlined, e.g., in the monographs \cite{AFP:2000_1} or \cite{Gi:1984_1}.\\

It will turn out in Section \ref{approx} and in Section \ref{exist} that the functional $K$ together with the class $\cbv$ is
the suitable choice in the setting at hand: in Section \ref{gen} we show that there exists a solution $\ovu$ of problem \reff{probbv}. 
Moreover, with the help of the geometric approximation procedure of Section \ref{approx}, we show 
in Corollary \ref{gen cor 1} that the infima of \reff{probsob} and \reff{probbv}
are equal. In the case that the superlinear part is given by a $N$-function, we obtain in addition a complete dual point of view.\\

In Section \ref{apriori} the apriori higher integrability and regularity results of the recent paper \cite{BF:2020_3} on splitting-type variational 
problems with linear growth are essentially refined and carried over to the mixed linear- superlinear setting.\\

Concerning the function $f_1$ we suppose that there exist real numbers $\mu >1$ and $\gamma \geq 0$ such that
\begin{equation}\label{intro ell 1}
c_1 \big(1+|t|\big)^{-\mu} \leq f_1''(t) \leq \overline{c}_1 \big(1+|t|\big)^\gamma \psp ,\fsp t \in \rz\psp ,
\end{equation} 
holds with constants $c_1$, $\overline{c_2} >0$. We note that $\mu > 1$ is in accordance with the required linear growth
of $f_1$.\\

For $f_2$ we suppose that there exist real numbers $\muhat < 2$ and $q \geq 1$ such that
\begin{equation}\label{intro ell 2}
c_2 \big(1+|t|\big)^{-\muhat} \leq f''_2(t) \leq \overline{c}_2 \big(1+|t|\big)^{q-2}  \psp ,\fsp t\in \rz\psp ,
\end{equation}
holds with constants $c_2$, $\overline{c_2} >0$. Since $f_2$ is of superlinear growth, the condition $\muhat < 2$ is
a quite mild assumption. In the case of ($p$,$q$)-growth of $f_2$ we have $-\muhat =p-2$ for some $p>1$. 
We also note that for $n >2$ the condition \reff{intro ell 2} is replaced by considering a function $f_2$: $\rz^{n-1}\to \rz$ such that
with constants $\lambda$, $\Lambda >0$ and for all $\xi_2$, $\eta \in \rz^{n-1}$
\begin{equation}\label{intro not 7}
\lambda \big(1+|\xi_2|\big)^{-\muhat} |\eta|^2 \leq D^2f_2(\xi_2)(\eta,\eta) \leq \Lambda \big(1+|\xi_2|\big)^{q-2}|\eta|^2 \psp .
\end{equation}

A Caccioppoli-type inequality w.r.t.~$\partial_2 \nabla u$ and $\partial_1 \nabla u$, respectively,
using in addition negative exponents gives different variants of regularity results 
depending on the properties of $f_{1,2}$.\\

In Section \ref{uni} we finally turn our attention to the question of uniqueness of solutions. 
First results were already given in Theorem \ref{apri 1818 theo 1} and Corollary \ref{apri 1818 cor 1}
by quoting \cite{Bi:1818}.\\

It remains to discuss the $N$-function case. 
If the ellipticity parameter from \reff{intro ell 1} satisfies $\mu < 2$, then the smoothness properties of $\sigma$ together with the uniqueness of $\sigma$
imply the uniqueness of generalized solutions. In order to make this argument precise we prove a generalization of 
\cite{Bi:2000_1}, Theorem 7, to the situation at hand.\\

\section{Approximation procedure}\label{approx}

In this section we present an approximation procedure which is adapted to the particular linear- superlinear setting.
Although the arguments seem to be quite technical, the principle idea is a geometric one.\\

We have to take care of various aspects:
\begin{itemize}
\item A retracting and smoothing procedure of the form $u_0+\eta_\eps*\big[(u-u_0)(x+\delta e_2)\big]$ 
is compatible with Lebesgue spaces. However it does not
work w.r.t.~the ``$\bv$''-direction $e_1$ which is due to the possible concentration of masses on the boundary.

\item In the linear growth situation the methods of local approximation (compare, e.g., \cite{Gi:1984_1}, Theorem 1.17, p.~14) 
together with some extension by $u_0$ outside of $\Omega$ (see, e.g., \cite{GMS:1979_1}, \cite{GMS:1979_2}) serve as a powerful tool.
However, a partition of the unity $\{\varphi_i\}$ is involved in this kind of argument.
This causes serious difficulties proving the convergence
of $f_2(\partial_2 w_m)$ for the approximating sequence $w_m$ since the derivatives of $\varphi_i$ do not cancel when calculating
the integral of $f_2$ evaluated at the corresponding expression.

\item Combining and adjusting both methods and using the geometric structure of the domain we obtain a partition of the unity such that
the derivatives w.r.t.~the relevant direction vanish.  
\end{itemize}

We start with a generalization of Lemma B.1 of \cite{Bi:1818} including strong $L^p$-convergence of $\partial_2 w_m$.
The main new feature is the way of constructing the sequence $\{w_m\}$ which is crucial for proving
Lemma \ref{lemma newapprox 1}.

\begin{lemma}\label{lemma b 1}
Let $w\in \bv(\Omega)$ such that $\partial_2 w \in L^p(\Omega)$ for some $1\leq p < \infty$ and such that 
$(w-u_0)\nu_2=0$ a.e.~on $\partial \Omega$.\\ 

Then there exits a sequence $\{w_m\}$ such that for all $m\in \nz$ we have
$w_m \in W^{1}_1(\Omega) \cap C^\infty (\Omega)$, $\partial_2 w_m \in L^p(\Omega)$, $\op{trace} w_m = \op{trace} w$ and
such that we the convergences
\begin{eqnarray*}
\lim_{m\to \infty} \iom |w_m -w|\dx &=& 0 \psp ,\\[2ex]
\lim_{m\to \infty} \iom \sqrt{1+|\nabla w_m|^2} \dx &=& \iom \sqrt{1+|\nabla w|^2} \psp ,\\[2ex]
\lim_{m\to \infty} \iom |\partial_2 w_m - \partial_2 w|^p \dx &=& 0 \psp .
\end{eqnarray*}
\end{lemma}

\emph{Proof.} Recalling our assumption \reff{intro not normal} imposed on the domain $\Omega$ we may consider w.l.o.g.~the case
\begin{equation}\label{b 1 1}
\Omega = (-1,1) \times (-1,1) \psp .  
\end{equation}

We fix a function $w\in BV(\Omega)$ and proceed in five steps.\\

{\bf Step 1.} In the following we suppose that $u_0=0$. The general case is obtained by considering $w-u_0$ and adding $u_0$ 
at the end of the proof.\\

We then reduce the problem by choosing
two smooth functions $\psi_1$, $\psi_2$: $[-1,1] \to [0,1]$ such that $\psi_1+\psi_2 \equiv 1$, 
$\psi_1(t) =0$ on $[-1,-1/2]$,  $\psi_1 (t) =1$ on $[1/2,1]$ 
and $\psi_2 (t) =1$ on $[-1,-1/2]$, $\psi_2(t) =0$ on $[1/2,1]$.\\

We consider $\psi_1(x_2) w$ and $\psi_2(x_2) w$ separately, hence w.l.o.g.~$w\equiv 0$  
in a neighborhood of $[x_2=-1]$.\\

{\bf Step 2.} Fix some $\eps_0 >0$ and let (w.r.t.~the $x_2$-direction)
\begin{equation}\label{b 1 2}
w_{\eps_0}(x) = w(x+\eps_0 e_2) \psp ,
\end{equation}
where $w$ is extended by $0$ on $(-1,1)\times [1,\infty)$. At the end of our proof we pass to the limit $\eps_0 \to 0$.\\

Thus we may suppose w.l.o.g. that
\begin{equation}\label{b 1 3}
w \equiv 0 \fsp\mbox{on}\fsp \Big[(-1,1) \times (1-\eps_0,1)\big] \cap 
\big[(-1,1)\times (-1,-1+\eps_0)\Big] \psp . 
\end{equation}

\vspace*{2ex}
{\bf Step 3.} We now take \cite{Gi:1984_1}, proof of Theorem 1.17, as a reference (compare also \cite{Bi:1818}),
Lemma B.1), fix $\eps >0$, (recalling Step 1 and Step 2) and for $l \in \nz$ we let
\[
\Omega_k = \Omega_k^l := \Bigg\{x\in \Omega:\psp -1 + \frac{1}{l+k} < x_1 < 1 - \frac{1}{l+k}\Bigg\} \psp , \msp k \in \nz_0\psp ,
\]
where $l$ is chosen sufficiently large such that
\begin{equation}\label{b 1 4}
\int_{\Omega-\Omega_0} |\nabla w|\dx < \eps \psp .
\end{equation}
With this notation we define $A_1 := \Omega_2$ and 
\begin{eqnarray*}
A_i = \Omega_{i+1}-\overline{\Omega}_{i-1} &:=& \Bigg\{x\in \Omega:\psp -1 + \frac{1}{l+i+1} < x_1 < -1 + \frac{1}{l+i-1}\\[2ex]
&&\qquad \mbox{and}\msp 1- \frac{1}{l+i-1} < x_1 <1 - \frac{1}{l+i+1}\Bigg\}\\[2ex]
&=:&\big\{x\in \Omega: \psp x_1 \in I_{i}^{-}\cup I_{i}^{+}\big\}\psp .
\end{eqnarray*}
A partition $\{\varphi_i\}$ of the unity is defined w.r.t.~these sets by
\[
\varphi_i \in C_0^\infty(A_i)\psp ,\fsp 0 \leq \varphi_i \leq 1 \psp , \fsp \sum_{i=1}^\infty \varphi_i = 1 \msp\mbox{on}\msp \Omega\psp . 
\]

For proving Lemma \ref{lemma newapprox 1} below it will be crucial to observe that the functions $\varphi_i$ 
may be chosen respecting the structure of the stripes, i.e.~for all $i\in \nz$
\begin{equation}\label{b 1 5}
\varphi_i (x_1,x_2) = \tilde{\varphi}_i (x_1)\psp , \fsp \tilde{\varphi}_{i}\in C^\infty_0\big(I_{i}^{-}\cup I_{i}^{+}\big) \psp .
\end{equation}

\vspace*{2ex}
{\bf Step 4.} Now we  proceed essentially as described in Lemma B.1 of \cite{Bi:1818}: 
let $\Omega_{-1}= \emptyset$ and denote by $\eta$ a smoothing kernel. On account of \reff{b 1 3} we select $\eps_i$ small enough
such that the smoothing procedure is well defined and such that we have
\begin{eqnarray}\label{b 1 6}
\op{spt} \eta_{{\eps_i}}* (\varphi_iw) &\subset& \Omega_{i+2}-\overline{\Omega}_{i-2}\psp ,\nonumber\\[2ex]
\iom \big|\eta_{{\eps_i}}*(\varphi_iw) - \varphi_i w\big| \dx &<& 2^{-i} \eps \psp ,\nonumber\\[2ex]
\iom \big|\eta_{{\eps_i}}* (w \nabla \varphi_i) - w \nabla \varphi_i\big| \dx&<& 2^{-i} \eps \psp ,\nonumber\\[2ex]
\iom \big|\eta_{{\eps_i}} *\partial_2 (\varphi_i w) - \partial_2 (\varphi_i w)\big|^p\dx &<& 2^{-i} \eps \psp . 
\end{eqnarray}
Moreover, the analogue to \reff{b 1 3} holds for $w_m$ with some $\tilde{\eps}_0 < \eps_0$. Here with the choice $\eps =1/m$ we have set
\[
w_m = \sum_{i=1}^{\infty} \eta_{{\eps_i}} * (\varphi_i w) \psp .
\]
By the above remarks we suppose with a slight abuse of notation (relabeling $\eps_0$) that we have in addition to
\reff{b 1 3} for all $m\in \nz$
\begin{equation}\label{b 1 7}
w_m \equiv 0 \fsp\mbox{on}\fsp \Big[(-1,1) \times (1-\eps_0,1)\big] \cap 
\big[(-1,1)\times (-1,-1+\eps_0)\Big] \psp . 
\end{equation}

Given \reff{b 1 7} we follow exactly the proof of Lemma B.1, where in particular the notion of a convex function $g$ of a measure
(see \cite{DT:1984_1}) is exploited via the representation
\[
\int_U g(\nabla w) := \sup_{\varkappa \in C^\infty_0(U:\rz^n), \psp |\varkappa| \leq 1}
\Bigg\{- \int_U w\op{div}\varkappa \dx - \int_U g^*(\varkappa) \dx \Bigg\} 
\]
and where $g$ is of linear growth and $g^*$ denotes the conjugate function (see the definition given in Section \ref{dual}).\\

{\bf Step 5.} With $\eps \ll \eps_0$ (i.e.~choosing $m=m(\eps_0)$ sufficiently large) we pass to the limit $\eps_0 \to 0$, 
which finally proves the lemma. \qed\\

Following the lines of Lemma \ref{lemma b 1} we obtain the convergence of the superlinear part of
the energy under consideration.

\begin{lemma}\label{lemma newapprox 1}
Given the notation of Lemma \ref{lemma b 1} we suppose that we have \reff{intro not 2} and \reff{A growth}.
Moreover, we now assume \reff{intro not 5}.\\ 

Then the sequence $\{w_m\}$ of Lemma \ref{lemma b 1} satisfies
\[
\iom f_2 (\partial_2 w_m) \dx \to \iom f_2(\partial w) \dx\fsp\mbox{as}\msp m\to \infty \psp .
\]
\end{lemma}

\emph{Proof.} We start with the first three steps of the proof of Lemma \ref{lemma b 1}, in particular we
have \reff{b 1 3}, \reff{b 1 7} and \reff{b 1 5}.\\

If $p$ is the exponent given in \reff{A growth}, then the strong $L^p$-convergence of the sequence $\{\partial_2 w_m\}$ 
yields (after passing to a subsequence)
\begin{equation}\label{newapprox 1 1}
\partial_2 w_m \to \partial_2 w \fsp\mbox{a.e.~in $\Omega$.}
\end{equation}

The first ingredient of the proof  follows from our assumption \reff{intro not 5} and Jensen's inequality, where we recall 
that in fact only finite sums are considered: for all $x\in \Omega$ and for all $m\in \nz$ we have
\begin{eqnarray}\label{newapprox 1 2}
f_2 \big(\partial_2 w_m\big) &=& f_2\Bigg(\partial_2 \sum_{i=1}^\infty \eta_{{\eps_i}}*(\varphi_i w)\Bigg)\nonumber
=  f_2\Bigg(\sum_{i=1}^\infty \eta_{{\eps_i}}*\partial_2(\varphi_i w)\Bigg)\nonumber\\[2ex]
&\leq & c  \sum_{i=1}^\infty  f_2\Big(\eta_{{\eps_i}}*\partial_2(\varphi_i w)\Big)
\leq  c  \sum_{i=1}^\infty \eta_{{\eps_i}}* f_2\big(\partial_2(\varphi_i w)\big)\psp .
\end{eqnarray}

We also recall \reff{b 1 5} which means $\partial_2 \varphi_i = 0$. In conclusion, \reff{newapprox 1 2} shows
\begin{equation}\label{newapprox 1 3}
f_2 (\partial_2 w_m) \leq c \sum_{i=1}^\infty \eta_{{\eps_i}}*f_2\big(\varphi_i \partial_2 w) \psp .
\end{equation}

Now we benefit from \reff{newapprox 1 1} and Egoroff's theorem: for any $\bar{\eps} >0$ 
and for any $i \in \nz$  there exists a measurable set $\aep$
such that 
\begin{equation}\label{newapprox 1 4}
|A_i - \aep| < \bar{\eps}_i \ll \bar{\eps}\fsp\mbox{and}\fsp
\partial_2 w_m \rightrightarrows \partial_2 w\fsp\mbox{on}\msp \aep \psp .
\end{equation}
A suitable choice of $\bar{\eps}_i$ is made in \reff{newapprox 1 6}. \\

With the help of  \reff{newapprox 1 3} one obtains for fixed $i\in \nz$ 
(note that by the first condition of \reff{b 1 6}, there exist at most three different numbers $k\in \nz$ such that the 
function $\eta_{\eps_k}*f_2(\varphi_k\partial_2w)  \not\equiv 0$  on $A_i$)
\begin{eqnarray}\label{newapprox 1 5}
\lefteqn{\int_{A_i - \aep} f_2(\partial_2 w_m) \dx}\nonumber\\[2ex] 
& \leq & c \sum_{k=1}^{\infty} \int_{A_i - \aep} \eta_{{\eps_k}}*f_2(\varphi_k \partial_2 w)\dx \nonumber\\[2ex]
&=& c \sum_{k=1}^\infty  \int_{A_i -\aep} \Bigg[\int \eta_{{\eps_k}}(x-y) f_2 \Big( \big(\varphi_k \partial_2 w\big)(y)\Big)\dy \Bigg] \dx\nonumber\\[2ex]
&=&  c \sum_{k=1}^\infty  \int_{A_i -\aep} \Bigg[\int_{B_1} \eta(z) f_2 \Big( \big(\varphi_k\partial_2 w\big)(x-\eps_k z)\Big)\dz \Bigg] \dx\nonumber\\[2ex]
&\leq &   c \sum_{k=1}^\infty  \int_{B_1} \eta(z) \Bigg[\int_{T^k_{i,\bar{\eps}}} f_2 \Big( \big(\partial_2 w\big)(y)\Big)\dy \Bigg] \dz\psp ,
\end{eqnarray}
where it is abbreviated ($|U|$ denoting the Lebesgue measure of $U\subset \Omega$)
\[
T^k_{i,\bar{\eps}} := \big\{y = x- \eps_k z:\psp x \in A_i-\aep\big\}\psp , \msp\mbox{in particular}\msp  |T^k_{i,\bar{\eps}}| = |A_i -\aep|\psp .
\]

Now, since for fixed $i$ the sum is just taken over three indices, we may choose $\bar{\eps}_i$ sufficiently small and finally obtain from \reff{newapprox 1 5} (recalling $\iom f_2(\partial_2 w)\dx < \infty$)
\begin{equation}\label{newapprox 1 6}
\int_{A_i - \aep} f_2(\partial_2 w_m) \dx \leq 2^{-i} \bar{\eps}\psp .
\end{equation}

Decreasing $\bar{\eps}_i$, if necessary, it may also be assumed that
\begin{equation}\label{newapprox 1 7}
\int_{A_i - \aep} f_2(\partial_2 w) \dx \leq 2^{-i} \bar{\eps}\psp .
\end{equation}

By \reff{b 1 3} and \reff{b 1 7} we note once more that only finite sums have to be considered and 
recalling \reff{newapprox 1 4}, \reff{newapprox 1 6} and \reff{newapprox 1 7} we obtain
\begin{eqnarray*}
\lefteqn{\Bigg| \iom f_2(\partial_2 w_m) \dx - \iom f_2(\partial_2 w)\dx\Bigg|}\\[2ex]  
&&\leq \sum_{i=1}^{N_0} \int_{A_i}\Big|f _2(\partial_2 w_m) -f_2(\partial_2 w)\Big|\dx\\[2ex]
& &\leq \sum_{i=0}^{N_0} \int_{\aep} \Big|f_2(\partial_2 w_m) -f_2(\partial_2 w)\Big| \dx + 2\bar{\eps}\\[2ex]
&&\leq \sum_{i=0}^{N_0} \sup_{\aep}\Big| f_2(\partial_2 w_m)-f_2(\partial_2 w)\Big| \psp  |A_i| +2  \bar{\eps} 
\leq 3 \bar{\eps}
\end{eqnarray*}
provided that $m >  m_0$ with $m_0$ sufficiently large. This finishes the proof of the lemma. \qed\\ 

\begin{remark}\label{remark newapprox 1}
Now we can shortly discuss one main difference to the model of Hencky plasticty investigated in \cite{Te:2018_1}. There an approximation lemma
is formulated as Theorem 5.3 in Chapter II. The convergence of the deviatoric part in terms of 
the density $f$ considered there corresponds to the convergence of the square root 
in Lemma \ref{lemma b 1} which follows from the linear growth of $f$ and the notion of a convex function of a measure.\\

Our main difficulty is proving the convergence of the $f_2$-energy. In the case of the Hencky plasticity the analogue is just a consequence
of considering the intermediate topology (defined in formula (3.37), Chapter II, of \cite{Te:2018_1}) 
which respects the linear operator $\op{div}v$, see also Theorem 3.4 and formula (5.53), Chapter II, of \cite{Te:2018_1}. 
\end{remark}

Now we define
\[
\hat{\Omega} := (-2,2)\times (-1,1)
\]
and for $w\in \bv(\Omega)$ we let
\[
\hatw := \left\{\begin{array}{rcl}
w&\mbox{on}&\Omega \psp ,\\
u_0&\mbox{on}&\hat{\Omega}-\Omega \psp ,
\end{array}\right.
\]
where $u_0$ represents a Lipschitz extension of
our fixed boundary datum from the space $W^{1}_\infty(\Omega)$.\\

We then have the validity of an approximation result corresponding to Lemma B.2 of \cite{Bi:1818}. It can be seen as a kind of generalization
of Lemma \ref{lemma b 1} and Lemma \ref{lemma newapprox 1} where now $C^\infty$ is replaced by $C^\infty_0$.

\begin{lemma}\label{lemma newapprox 2}
Using the above notation suppose that $w\in\bv(\Omega)$, $\|\partial_2 w\|_{L_A(\Omega)} < \infty$ and that
we have \reff{intro not 1}, \reff{intro not 2}, \reff{A growth} and \reff{intro not 5}. Then there exists a sequence $\{w_m\}$ in $u_0+C^\infty_0(\Omega)$ such that
passing to the limit $m\to \infty$ we have
\[
\begin{array}{cc}
i)&\dis \hatw_m \to \hat{w} \fsp\mbox{in}\msp L^1\big(\hat{\Omega}) \psp ,\\[2ex]
ii)&\dis \int_{\hat{\Omega}} \sqrt{1+|\nabla \hatw_m|^2}\dx \to \int_{\hat{\Omega}} \sqrt{1+|\nabla \hatw|^2}\psp ,\\[3ex]
iii)&\dis \iom f_2(\partial_2 w_m)\dx \to \iom f_2(\partial_2 w) \dx \psp .
\end{array} 
\]
\end{lemma}

\section{Existence of solutions}\label{exist}

There are two approaches towards the existence of generalized solutions to problem \reff{probsob}. 
The first one follows the direct method and leads to the 
existence of solutions to problem \reff{probbv}. This works under quite weak assumptions,
for example, the densities $f_1$ and $f_2$ need not to be of class $C^2$.\\

The second approach yields the stress tensor as the unique solution of the dual problem and by the
stress- strain relation a complete picture of the situation is drawn. 
However, following the duality approach, we have to suppose that $f_2$ is given in terms of a $N$-function.

\subsection{Generalized solutions}\label{gen}

\begin{theorem}\label{gen theo 1}
Suppose that we have \reff{intro not 1} - \reff{intro not 5}.  
Then the relaxed problem \reff{probbv} admits a solution $\ovu \in \cbv$.
\end{theorem}

\emph{Proof of Theorem \ref{gen theo 1}.} We recall that w.l.o.g.~we suppose $n=2$ and
consider a $K$-minimizing sequence $\big\{u\ob{n}\big\}$ in the
admissible  class $\cbv$ of comparison functions.
After passing to a subsequence we may assume that there exits a function $\ovu \in \bv(\Omega)$ and a function $v \in L^p(\Omega)$
such that as $n \to \infty$
\begin{equation}\label{gen 1}
u\ob{n} \to \ovu \msp\mbox{in}\msp L^1(\Omega)\psp ,\fsp
\partial_2 u\ob{n} \rightharpoondown v \msp\mbox{in}\msp L^p(\Omega) \psp .
\end{equation}
Here, in the case $p=1$, we refer to the Theorem of De LaValee-Poussin (see, e.g., \cite{KR:1961_1}).\\

We have for any $\varphi \in C_0^\infty(\Omega)$
\[
\iom u\ob{n} \partial_2 \varphi \dx = - \iom \partial_2 u\ob{n} \varphi \dx \psp ,
\]
hence $v = \partial_2 \ovu$ and since we have for any $\psi \in C^\infty(\Omega)$
\begin{eqnarray*}
\iom u\ob{n} \partial_2 \psi \dx &=& - \iom \partial_2 u\ob{n} \psi \dx + \int_{\partial \Omega} u_0 \psi \nu_2 \D\mathcal{H}^1\psp ,\\[2ex]
\iom \ovu \partial_2 \psi \dx &=& - \iom \partial_2 \ovu  \psi \dx + \int_{\partial \Omega} \ovu \psi \nu_2 \D\mathcal{H}^1\psp ,
\end{eqnarray*}
the convergences stated in \reff{gen 1} prove $\ovu \in \cbv$.\\

We note that
\[
\liminf_{n\to \infty} K[u\ob{n}] \geq \liminf_{n\to \infty} K_1\big[u\ob{n}\big] + \liminf_{n\to \infty}E\big[\partial_2 u\ob{n}\big] \psp .
\]
By \cite{Re:1968_1}, see also \cite{AFP:2000_1}, Theorem 5.47, p.~304, we have the lower semicontinuity
\[
K_1\big[\ovu\big] \leq \liminf_{n\to \infty} K_1\big[u\ob{n}\big] \psp .
\]
Discussing $E$ we cite Theorem 2.3, p.~18, of \cite{Gia:1983_1}, hence
\[
E\big[\partial_2\ovu\big] \leq \liminf_{n\to\infty} E\big[\partial_2 u\ob{n}\big] \psp .
\]
Since $\big\{u\ob{n}\big\}$ was chosen as a $K$-minimizing sequence, the proof of Theorem \ref{gen theo 1} is complete. \qed\\

Now, on account of our approximation Lemma \ref{lemma newapprox 2}, we have 

\begin{cor}\label{gen cor 1}
With the notation and under the hypotheses of Theorem \ref{gen theo 1} we have
\[
\inf_{w\in \csob}J[w] = \inf_{v\in \cbv}K[v] = K[\ovu]\psp .
\]
\end{cor}

\subsection{The dual solution}\label{dual}

Another approach leading to an analogue of the stress tensor, occurring as basic quantity in problems from mechanics,
is to consider the dual problem. As the main references on convex analysis 
we mention \cite{Ro:1970_1} and \cite{ET:1976_1}.\\ 

Let us assume that we have \reff{intro not 2} with $A\big(|\xi_2|\big) = f_2(\xi_2)$  
for all $\xi\in \rz^{n-1}$, with $A$ being of class $C^1\big([0,\infty)\big)$ and with $A$ satisfying $(N1)$-$(N3)$.
In this case we suppose for notational simplicity that $f$: $\rz^2 \to \rz$,
\begin{equation}\label{dual 1}
f(\xi) = f_1(\xi_1) + A\big(|\xi_2|\big) \psp , \fsp \xi \in \rz^2\psp .
\end{equation}

As usual we define the conjugate function $A^*$: $[0,\infty) \to [0,\infty)$ by
\begin{equation}\label{dual 2}
A^{*}(s) := \max_{t \geq 0} \big\{st - A(t)\big\} \psp .
\end{equation}
and note that we have for all $t \in [0,\infty)$
\begin{equation}\label{dual 3}
A(t)+ A^{*}\big(A'(t)) = t A'(t) \psp .
\end{equation}
In order to obtain a well-posed dual problem we additionally require
\begin{equation}\label{dual 4}
A^{*} \big(A'(t)\big) \leq c \Big[A(t)+ 1 \Big]\fsp\mbox{for all}\msp t \in \rz \psp .
\end{equation}

Since
\[
f_1^{*}(s) := \sup_{\xi_1 \in \rz} \Big\{s\xi_1 - f_1(\xi_1)\Big\} \psp ,
\]
we obtain from the decomposition of $f$ 
\begin{equation}\label{dual 5}
f^{*}(\xi) =f_1^*(\xi_1) + A^*(|\xi_2|) 
\end{equation}
as formula for the conjugate function $f^*$: $\rz^2 \to \rz$.
The conjugate function $f^*$ satisfies in correspondence to \reff{dual 3}
\begin{equation}\label{dual 6}
f(\xi) + f^*\big(Df(\xi)\big) = \xi \cdot Df(\xi) \psp , \fsp \xi \in \rz^2\psp .
\end{equation}
Given these preliminaries we define the Lagrangian
\begin{eqnarray}\label{dual 7}
l(v,\tau) &:=& \iom \tau \cdot \nabla v \dx - \iom f^{*}_1(\tau_1) \dx -   \iom A^*\big(|\tau_2|\big)\dx \psp ,\nonumber\\[2ex]
&&v \in u_0+ \wnullaa\psp ,\fsp \tau\in L^{\infty,A^*}(\Omega;\rz^2) \psp .
\end{eqnarray}
In \reff{dual 7} we have set
\[
L^{\infty,A^*} (\Omega;\rz^2) := L^\infty(\Omega)\times L_{A^*}(\Omega) \psp .
\]
With the help of the formula for the conjugate function given in \reff{dual 5} we have the representation for the energy $J$
defined in \reff{probsob}
\begin{eqnarray}\label{dual 8}
J[w] &=&\sup_{\varkappa \in L^{\infty,A^*}(\Omega;\rz^2)}\Bigg\{\iom \varkappa \cdot \nabla w \dx  - \iom f_1^*(\kappa_1) \dx\nonumber\\[2ex]
&&\hspace*{4cm} -\iom A^*\big(|\kappa_2|\big) \dx \Bigg\}\psp , \nonumber\\[2ex]
&=&\sup_{\varkappa \in L^{\infty,A^*}(\Omega;\rz^2)}l(w,\varkappa)\psp , \fsp w \in u_0+\wnullaa\psp .
\end{eqnarray}
The dual functional finally is defined via
\begin{equation}\label{dual 9}
R[\tau] := \inf_{w\in u_0 + \wnullaa}l(w,\tau)\psp ,\fsp \tau\in L^{\infty,A^*}(\Omega;\rz^2) \psp .
\end{equation}
This functional leads to the dual problem as the maximizing problem
\begin{equation}\label{dual 10}
R[\tau] \to \max \fsp\mbox{in}\msp \tau \in L^{\infty,A^*}(\Omega;\rz^2)\psp .
\end{equation}

Then we have recalling Theorem \ref{gen theo 1}

\begin{theorem}\label{dual theo 1}
Suppose that we have our general assumptions \reff{intro not 1} -  \reff{intro not 5}. 
Moreover, suppose that $f$ is given in \reff{dual 1} with $A$ satisfying \reff{dual 4}.
Let $\ovu$ denote a solution of the problem \reff{probbv}.\\ 

Then the ``stress tensor'' defined by
\begin{equation}\label{dual 11}
\sigma(x) := D f\big(\nabla^a \ovu\big) = \Big(f_1'\big(\partial_1^a \ovu\big), A'\big(|\partial_2 \ovu|\big)\Big) 
\end{equation}
is of class $L^{\infty,A^*}(\Omega;\rz^2)$ and maximizes the dual variational problem \reff{dual 10} with $R$ given in \reff{dual 9}. 
\end{theorem}

\emph{Proof.} We first note that the boundedness of $|f_1'|$ and condition \reff{dual 4} imply $\sigma \in L^{\infty,A^*}(\Omega;\rz^2)$.\\

We then follow an Ansatz similar to Lemma 5.1 of \cite{FT:2015_1}. For any 
$v\in u_0+\wnullaa$ we have recalling \reff{dual 7} und using \reff{dual 6}
\begin{eqnarray}\label{dual 12}
l(v,\sigma) &=& \iom \nabla v \cdot Df\big(\nabla^a \ovu\big) \dx - \iom f^{*}\Big(Df\big(\nabla^a \ovu\big)\Big) \dx\nonumber\\[2ex]
&=& \iom Df\big(\nabla^a \ovu\big) \cdot \big(\nabla v -\nabla^a \ovu\big)\dx + \iom f\big(\nabla^a \ovu\big) \dx \psp .
\end{eqnarray}

\vspace*{2ex}
Now given $|t| \ll 1$ let $\ovu_t := \ovu + t(v-\ovu) \in u_0+\wnullaa$. 
The $K$-minimality of $\ovu$ obviously implies
\[
\frac{\D}{\D t}_{|t=0}K[\ovu_t] = 0\psp ,
\]
hence by $\nabla^s v =0$
\begin{eqnarray}\label{dual 13}
0&=& \iom Df\big(\nabla^a \ovu\big)\cdot \big(\nabla v - \nabla^a \ovu) \dx 
+ \frac{\D}{\D t}_{|t=0}\iom f_1^\infty \Bigg(\frac{\partial_1^s \ovu_t}{\big|\partial_1^s \ovu_t\big|}\Bigg) \D \big|\partial_1^s \ovu_t\big|\nonumber \\[2ex] 
&&+\frac{\D}{\D t}_{|t=0}  \int_{\partial \Omega} f_1^\infty \Big( \big(u_0-\ovu_t\big) \nu_1\Big) \D \mathcal{H}^1 \psp ,
\end{eqnarray}
where $\partial_1^s \ovu_t = (1-t) \partial_1^s \ovu$. Now we note that
\[
\frac{\D}{\D t}_{|t=0} \iom f_1^\infty \Bigg(\frac{\partial_1^s \ovu}{\big| \partial_1^s \ovu\big|}\Bigg)\D \Big((1-t)\big|\partial_1^s \ovu\big|\Big) =
- \iom f_1^\infty \Bigg(\frac{\partial_1^s \ovu}{\big|\partial_1^s \ovu\big|}\Bigg) \D\big|\partial_1^s u\big| \psp ,
\]
and since $v$ takes the boundary data $u_0$ on $\partial \Omega$ we have
\[
\frac{\D}{\D t}_{|t=0} \int_{\partial \Omega} f_1^\infty \big( (u_0 -\ovu_t)\nu_1\big) \D\mathcal{H}^1
= - \int_{\partial \Omega} f_1^\infty \big( (u_0 -\ovu) \nu_1\big) \D \mathcal{H}^1 \psp .
\]
Hence, inserting \reff{dual 12} in \reff{dual 13} we have shown
\[
l(v,\sigma) = K[\ovu] \fsp\mbox{for any}\msp v \in u_0+\wnullaa
\]
and taking the infimum w.r.t.~the comparison function $v$ we have
\begin{equation}\label{dual 14}
R[\sigma] \geq K[\ovu]\psp .
\end{equation}
We already know from Section \ref{approx} that $\inf J = K[\ovu]$ and the representation \reff{dual 8}
finally yields ($w\in u_0 + \wnullaa$)
\begin{eqnarray*}
J[w]  &=& \sup_{\varkappa \in L^{\infty,A^*}(\Omega;\rz^2)} l(w,\varkappa)\\[2ex] 
&\geq& \sup_{\varkappa \in L^{\infty,A^*}(\Omega;\rz^2)}\Bigg\{\inf_{v \in u_0+\wnullaa} l(v,\varkappa)\Bigg\}\\[2ex]
&=& \sup_{\varkappa \in L^{\infty,A^*}(\Omega;\rz^2)} R[\varkappa]\psp ,\fsp\mbox{i.e.}\\[2ex]
\inf_{w \in u_0+\wnullaa}J[w] &\geq& \sup_{\varkappa \in L^{\infty,A^*}(\Omega;\rz^2)} R[\varkappa] \psp .
\end{eqnarray*}
This together with \reff{dual 14} and Corollary \ref{gen cor 1} proves the theorem. \qed\\

\section{Higher integrability}\label{apriori}

\subsection{Anisotropic behaviour of the superlinear part with $q=2$}\label{apriori gen}

Let us start by recalling Theorem 6.5 of \cite{Bi:1818} together with our general assumption $\muhat < 2$ in \reff{intro ell 2}. 

\begin{theorem}\label{apri 1818 theo 1}
Suppose that we are given the assumptions \reff{intro not 1}, \reff{intro not 2}, \reff{intro ell 1} with $\mu < 2$,
$\gamma =0$ and \reff{intro ell 2} with $q=2$. Let
\begin{eqnarray*}
u^* \in \mathcal{M}&:=&\Big\{ u\in \bv (\Omega):\psp \mbox{\rm $u$ is the $L^1$-limit of a $J$-minimizing}\\
&&\hspace*{3cm}\mbox{\rm sequence from $u_0 +\wnull$} \Big\}\psp .
\end{eqnarray*}
Then $u^*$ is of class $C^{1,\alpha}(\Omega)$ for any $0 < \alpha < 1$. Moreover, the elements of $\mathcal{M}$ are uniquely 
determined up to constants.
\end{theorem}

\begin{remark}\label{apri 1818 rem 1}
We again emphasize that the theorem holds without the restriction
$n=2$, i.e.~the general case $f_2$: $\rz^{n-1}\to \rz$
\[
\lambda \big(1+|\xi_2|\big)^{-\muhat} |\eta|^2 \leq D^2f_2(\xi_2)(\eta,\eta) \leq \Lambda |\eta|^2 \psp ,
\fsp \xi_2 ,\psp \eta \in \rz^{n-1}\psp ,
\]
$\lambda$, $\Lambda >0$ is included without being explicitely mentioned.
\end{remark}

\emph{Main idea of the proof of Theorem \ref{apri 1818 theo 1}.} The proof of the theorem is based on uniform apriori estimates 
for the minimizers $\udel$ of the standard quadratic regularization. However,
is not immediate that $\{\udel\}$ is a $J$-minimizing sequence, if the superlinear part is not generated by a $N$-function.\\

One way to overcome this difficulty is to introduce some kind of local regularized stresss tensor $\sigma_\delta$ without knowing that
a solution to a dual problem exists (see Remark 6.15 of \cite{Bi:1818}).\\ 

This, together with the equation $\op{div} \sigma_\delta =0$,
leads to the minimality of weak $L^1$-cluster points of $\{\udel\}$ via a generalized $\inf -\sup$ relation. In conclusion, a variational inequality is derived
for any $u^* \in \mathcal{M}$ which provides the regularity of $u^*$ by Corollary 6.13 of \cite{Bi:1818} in the $C^{0,\alpha}$-regularity
of the stress tensor. \hspace*{\fill}\qed\\

With our approximation result Corollary \ref{gen cor 1} we identify the elements of $\mathcal{M}$ with $K$-minimizers and observe that
the class $\cbv$ is defined respecting the condition $(w-u_0)\nu_2=0$ $\mathcal{H}^1$ a.e.~on $\partial \Omega$.

\begin{cor}\label{apri 1818 cor 1}
Given the assumptions of Theorem \ref{apri 1818 theo 1}, the problem \reff{probbv} is uniquely solvable.
\end{cor}

With the help of Theorem \ref{apri 1818 theo 1} we now establish the first regularity result of this paper which is very much in
the spirit of \cite{BF:2020_3}, i.e.: if we drop the ellipticty condition $\mu < 2$, then we still have higher integrability of
$\partial_2 u$.

\begin{theorem}\label{apri 1818 theo 2}
Suppose that we have the assumptions of Theorem \ref{apri 1818 theo 1} without
the requirement $\mu < 2$.\\

Then there exists a generalized minimizer $u \in \mathcal{M}$ such that
\[
\partial_2 u \in L^\chi_{\op{loc}}(\Omega)\fsp\mbox{for any finite $\chi$}\psp .
\]
\end{theorem}

\emph{Proof of Theorem \ref{apri 1818 theo 2}.} With the ideas presented in \cite{BF:2020_3} in the linear growth case,
the main point is to introduce in the mixed linear-superlinear case a suitable regularization procedure
to obtain a sufficiently smooth minimizing sequence.\\

We fix $1 < \nu < 2$ and  define (compare \cite{BF:2003_1} and related papers)
\begin{eqnarray*}
\Phi_\nu(t) &:=& (\nu-1) \int_0^t\int_0^s (1+r)^{-\nu}\dr\ds\psp , \\[2ex]
 &=&  t - \frac{1}{2-\nu} (1 + t)^{2- \nu} -  \frac{1}{\nu - 2} \psp .
\end{eqnarray*}
Then $\Phi_\nu$ satisfies \reff{intro ell 1} with $\mu = \nu$ and $\gamma = 0$ (in fact we may even choose $\gamma =-1$).\\

With
\[
f_\delta(\xi) = \delta \Phi_\nu(|\xi_1|)  + f_1(\xi_1) + f_2(\xi_2)
\]
we consider with the obvious meaning of notation the regularized minimization problem ($\delta \ll 1$)
\begin{equation}\label{apri 1818 1}
K_\delta[w] = K_{1,\delta}[w] + E\big[\partial_2 w\big]\to\min \fsp\mbox{in the class}\msp \cbv\psp .
\end{equation}
By Corollary \ref{apri 1818 cor 1} there exists a unique solution denoted by $\udel$ and by Theorem \ref{apri 1818 theo 1}
$\udel$ is of class $C^{1,\alpha}(\Omega)$, hence arguing with the Euler equation we also have
$\udel \in W^{2}_{2,\op{loc}}(\Omega)$. \\

For the minimizing property of the sequence $\{\udel\}$ we observe for any fixed $w \in \csob$ 
\[
J[\udel] \leq J_\delta[\udel] \leq J_\delta[w]\psp ,
\]
hence 
\[
\liminf_{\delta \to 0} J[\udel] \leq \limsup_{\delta \to 0} J_\delta[w] = J[w]\psp .
\]
Passing to a subsequence, if necessary,  $\{\udel\}$ obviously is a $J$-minimizing sequence.\\

Now we argue exactly as in \cite{BF:2020_3}, proof of Theorem 1, and complete the proof of Theorem \ref{apri 1818 theo 2}.\hspace*{\fill}\qed\\

\subsection{Higher integrability results in case of superlinear parts of isotropic $p$-growth}\label{apriori p}

In this subsection we restrict our considerations to the 
case that the superlinear part is of $p$-power growth. This is done in order to simplify our exposition which otherwise
would be based on additional parameters. Generalizations are left to the reader.\\ 

The main issue is that linear growth conditions are also compatible with $p$-growth conditions, $p >2$.
Higher integrability results for $\partial_2 u$ are established in Theorem \ref{apri theo 1} under quite general assumptions.\\

In Theorem \ref{apri theo 2} we suppose in addition that $\mu < 2$ and obtain higher integrability of the full gradient,
again withount restrictions on $p$.\\

Throughout this section we concentrate on the following model case:\\

Given the general splitting hypothesis \reff{intro not 1}, we suppose that $f_1$ is of linear growth,
\[
a_1 |t| - a_2\leq  f_1(t)\leq  a_3 |t| + a_4\psp , \fsp t \in \rz\psp ,
\]
$a_2$, $a_4 \geq 0$, $a_1$, $a_3>0$, satisfying for some $\mu >1$, $\gamma \geq 0$
\begin{equation}\label{apri 1}
c_1 \big(1+|t|\big)^{-\mu} \leq f_1''(t) \leq \bar{c}_1 \big(1+|t|\big)^\gamma \psp , \fsp t \in \rz \psp ,
\end{equation}
with constants $c_1$, $\bar{c}_1 > 0$.\\

We recall the examples discussed in \cite{BF:2020_3} having linear growth and with unbounded second derivative,
i.e.~the case $\gamma >0$ is included in our considerations.\\

The function $f_2$ is supposed to be of $p$-power growth, $p>1$, in the sense that 
\[
b_1 |t|^p - b_2\leq  f_2(t)\leq  b_3 |t|^p + b_4\psp , \fsp t \in \rz\psp ,
\]
$b_2$, $b_4 \geq 0$, $b_1$, $b_3>0$, satisfying in addition
\begin{equation}\label{apri 2}
c_2 \big(1+|t|\big)^{p-2} \leq f_2''(t) \leq \bar{c}_2 \big(1+|t|\big)^{p-2}\psp , \fsp t \in \rz \psp ,
\end{equation}
with constants $c_2$, $\bar{c}_2 > 0$.\\

\begin{theorem}\label{apri theo 1}
Suppose that the hypotheses of this section are valid.
\begin{enumerate}
\item If $\gamma =0$, then there exists a generalized minimizer $u\in \mathcal{M}$ such that
\[
\partial_2 u \in L^\chi_{\op{loc}}(\Omega) \fsp\mbox{for all finite}\msp \chi > p\psp .
\]
\item If 
\begin{equation}\label{apri 3}
0 \leq \gamma < \frac{p}{p+1} \psp ,
\end{equation}
then there exists a generalized minimizer $u\in \mathcal{M}$ such that
\[
\partial_2 u \in L^\chi_{\op{loc}} (\Omega)\fsp\mbox{for some}\msp \chi > p+1 \psp .
\]
\end{enumerate}
\end{theorem}

\emph{Proof of Theorem \ref{apri theo 1}.} The proof of Theorem \ref{apri theo 1} is carried out in four steps.\\

\emph{Step 1. Regularization.} We fix some $0 < \delta < 1$ and define an appropriate variant
of regularization: let
\[
f_\delta(\xi)  :=  \delta \big(1+|\xi_1|^2\big)^{\frac{p}{2}} + f_{1}(\xi_1) + f_{2}(\xi_2)
\psp ,  \fsp \xi =(\xi_1,\xi_2) \in \rz^2 \psp .
\]
We note that this Ansatz respects the splitting structure of the energy density under consideration.\\

We then consider the  minimization problem
\renewcommand{\theequation}{\mbox{1.1$_\delta$}}
\begin{equation}\label{delta}
J_\delta[w] := \iom f_\delta(\nabla w) \dx \to \min\fsp\mbox{in}\msp u_0 + W^{1,p}_{0}(\Omega)
\end{equation}\addtocounter{equation}{-1}\renewcommand{\theequation}{\mbox{\arabic{section}.\arabic{equation}}}%
with $u_\delta$ denoting the unique solution of \reff{delta} satisfying in addition (see the monographs 
\cite{Gi:2003_1}, Theorem 8.1, p.~267, and \cite{LU:1968_1}, Theorem 5.2, p.~277)
\begin{equation}\label{reg 0}
\udel \in W^{2,2}_{\op{loc}}(\Omega) \cap W^{1,\infty}_{\op{loc}}(\Omega) \psp .
\end{equation}
In the situation at hand one may use the dual problem in order to show that the sequence $\{\udel\}$
is $J$-minimizing. To this purpose we let
\begin{eqnarray*}
\tau_\delta &:=& \nabla f\big(\nabla u_\delta\big)\psp , \\
\sigma_\delta &:=& \delta X_\delta + \tau_\delta = \nabla f_\delta\big(\nabla \udel\big)\psp ,
\fsp X_\delta := p \big(1+|\partial_1 \udel|^2\big)^\frac{p-2}{2} \partial_1 \udel \psp ,
\end{eqnarray*}
and  adapt the arguments of \cite{Bi:1818}, Section 4.1.2: we note that $\sigma_\delta$ is of class $W^{1,2}_{\op{loc}}$ 
(recall \reff{reg 0}) satisfying $\op{div}\sigma_\delta =0$ and by 
$J_\delta[\udel] \leq J_\delta[u_0] \leq J_1[u_0]$ we have with a finite constant $c$ independent of $\delta$ and not relabelled in the following
\begin{eqnarray}\label{reg 1}
\delta \iom \big(1+|\partial_1 \udel|^2\big)^\frac{p}{2} \dx \leq c \psp , && \mbox{i.e.}\fsp 
\|\delta^{\frac{p-1}{p}} X_\delta\|_{p/(p-1)}\leq c\\
\iom f(\nabla \udel)\dx \leq c\psp ,&&
\|\tau_\delta\|_{p/(p-1)} \leq c\psp . \label{reg 2}
\end{eqnarray}
Observe that \reff{reg 1} implies as $\delta \to 0$
\begin{equation}\label{reg 3}
\delta X_\delta \rightharpoondown 0 \fsp\mbox{in}\msp L^{\frac{p}{p-1}}(\Omega) \psp .
\end{equation}
After passing to subsequences we obtain from \reff{reg 2} and \reff{reg 3} as $\delta \to 0$
\begin{equation}\label{reg 4}
\tau_\delta\psp , \msp \sigma_\delta \rightharpoondown: \sigma\fsp\mbox{in}\msp L^{\frac{p}{p-1}}(\Omega) \psp .
\end{equation} 

We recall the equation
\[
\tau_\delta : \nabla \udel - f^*(\tau_\delta) = f(\nabla \udel)
\]
and arrive at
\begin{eqnarray*}
J_\delta[\udel] &=& \delta \iom \big(1+|\partial_1 \udel|^2\big)^{\frac{p}{2}}\dx + \iom \big[\tau_\delta:\nabla \udel -f^*(\tau_\delta)\big]\dx\\
&=&  \delta \iom \big(1+|\partial_1 \udel|^2\big)^{\frac{p}{2}}\dx + \iom \big[\sigma_\delta:\nabla \udel -f^*(\tau_\delta)\big]\dx\\
&&-\delta p \iom \big(1+|\partial_1 \udel|^2\big)^{\frac{p-2}{2}}|\nabla \udel|^2\dx
\end{eqnarray*}
Using $\op{div} \sigma_\delta =0$ we obtain
\begin{eqnarray*}
J_\delta[\udel] &=&  \delta \iom \big(1+|\partial_1 \udel|^2\big)^{\frac{p}{2}}\dx + \iom \big[\sigma_\delta:\nabla u_0 -f^*(\tau_\delta)\big]\dx\nonumber\\
&&-\delta p \iom \big(1+|\partial_1 \udel|^2\big)^{\frac{p-2}{2}}|\nabla \udel|^2\dx \nonumber\\
&=& \iom \big[\tau_\delta:\nabla u_0 -f^*(\tau_\delta)\big]\dx + \delta \iom X_\delta : \nabla u_0 \dx\nonumber\\
&&+(1-p)\delta \iom \big(1+|\partial_1 \udel|^2\big)^{\frac{p}{2}}\dx
+ \delta p \iom \big(1+|\partial_1 \udel|^2\big)^{\frac{p-2}{2}}\dx \psp .
\end{eqnarray*}
Here, as $\delta \to 0$, the second integral (recall \reff{reg 3}) and the last integral  on the right-hand side
converge to $0$ and we obtain as $\delta \to 0$ (recall \reff{reg 4} and upper
semicontinuity of $-f^*$)
\begin{eqnarray*}
R[\varkappa] &\leq&  \inf_{u\in \csob} J[u] \leq J[\udel] \leq J_\delta [\udel]\\
&\to&  R[\sigma] + (1-p) \delta \iom \big(1+|\partial_1 \udel|^2\big)^{\frac{p}{2}}\dx \psp .
\end{eqnarray*}
Hence we have the minimizing property of the sequence $\{\udel\}$ and additionally
\[
\delta \iom \big(1+|\partial_1 \udel|^2\big)^{\frac{p}{2}} \dx \to 0\fsp\mbox{as}\msp \delta \to 0\psp .
\]

\emph{Step 2. Caccioppoli-type inequality.} Proceeding with the proof of Theorem \ref{apri theo 1} we note that by $u_0\in W^{1}_\infty(\Omega)$
\[
\sup_{\delta} \|\udel\|_{L^{\infty}(\Omega)} < \infty \psp .
\]
As usual we let
\[
\gam{i}{} : = 1+ |\partial_i\udel|^2\psp , \fsp i=1,\psp 2\psp .
\]

We note that for proving $ii$) of Theorem \ref{apri theo 1}, we also need an iterated Caccioppoli-type inequality
with negative exponents.

\begin{proposition}\label{apri prop 1}
Fix $l\in \nz$ and suppose that $\eta \in C^\infty_0(\Omega)$, $0 \leq \eta \leq 1$. Then, given the assumptions at the beginning of
Section \ref{apriori p}, the inequality
\begin{eqnarray}\label{apri 4}
\lefteqn{\iom D^2 f_\delta(\nabla \udel)\big(\nabla \partial_2 \udel, \nabla \partial_2 \udel\big) \eta^{2l} 
\gam{2}{\alpha} \dx}\nonumber \\ 
&& \leq c \iom D^2f_\delta(\nabla \udel) (\nabla \eta,\nabla \eta)\eta^{2l-2} \gam{2}{\alpha +1}\ \dx
\end{eqnarray}
holds for any  $\alpha > - 1/2$, which in particular implies (again for all $\alpha > -1/2$)
\begin{eqnarray}\label{apri 5}
\lefteqn{\iom \eta^{2l} \gam{2}{\alpha + \frac{p-2}{2}} |\partial_{22}\udel|^2 \dx}\nonumber\\
&\leq& c \Bigg[ 1+ \iom |\nabla \eta|^2 \eta^{2l-2}\gam{2}{\frac{\alpha+1}{1-\gamma}}\dx
+ \iom |\nabla \eta|^2 \eta^{2l-2}\gam{2}{\alpha +1+\frac{p-2}{2}} \dx \Bigg] \psp .
\end{eqnarray}
Here and in what follows $c$ is a finite constant independent of $\delta$ ($c=c(l,\alpha)$).
\end{proposition}

\emph{Proof of Proposition \ref{apri prop 1}.} Let us first suppose that $-1/2 < \alpha \leq 0$. \\

We differentiate the Euler equation
\[
0 = \iom Df_\delta (\nabla \udel) \cdot \nabla \varphi \dx \fsp\mbox{for all}\msp
\varphi\in C^\infty_0(\Omega)
\]
by inserting $\varphi = \partial_2 \psi$ as test function, hence
\[
0 = \iom D^2 f_\delta(\nabla \udel) \big(\nabla \partial_2 \udel, \nabla \psi\big) \dx \fsp
\mbox{for all}\fsp \psi \in C^\infty_0(\Omega) \psp .
\]
With the choice
\[
\psi := \eta^{2l} \partial_2 \udel \gam{2}{\alpha}
\]
we obtain
\begin{eqnarray}\label{apri 6}
\lefteqn{\iom D^2f_\delta(\nabla \udel) \big(\nabla \partial_2 \udel , \nabla \partial_2 \udel\big)
\eta^{2l} \gam{2}{\alpha}\dx}\nonumber\\
&=& - \iom D^2f_\delta(\nabla \udel)\big(\nabla \partial_2 \udel, \nabla \gam{2}{\alpha}\big)
\partial_2 \udel \eta^{2l}\dx\nonumber\\
&& - \iom D^2f_\delta(\nabla \udel)\big(\nabla \partial_2 \udel, \nabla (\eta^{2l})\big)
\partial_2 \udel \gam{2}{\alpha} \dx  =: S_1+S_2 \psp ,
\end{eqnarray}
where we note
\begin{eqnarray*}
|S_1| &=&  2 |\alpha| \iom D^2f_\delta (\nabla \udel) \big(\nabla \partial_2 \udel, \nabla \partial_2 \udel\big)
|\partial_2 \udel|^2\gam{2}{\alpha-1}\eta^{2l}\dx\\
&\leq & 2  |\alpha| \iom D^2f_\delta (\nabla \udel) \big(\nabla \partial_2 \udel, \nabla \partial_2 \udel\big)\gam{2}{\alpha}\eta^{2l}\dx \psp .
\end{eqnarray*}

Since $2 |\alpha| < 1$, we may absorb $|S_1|$ on the left-hand side of \reff{apri 6} with the result
\begin{equation}\label{apri 7}
\iom D^2f_\delta (\nabla \udel) \big(\nabla \partial_2 \udel, \nabla \partial_2 \udel\big) \eta^{2l}\gam{2}{\alpha}\dx\leq c |S_2| \psp .
\end{equation}

In the case $\alpha >0$ we immediately have \reff{apri 7} on account of the observation
\[
\iom D^2f_\delta(\nabla \udel)\big(\nabla \partial_2 \udel, \nabla \gam{2}{\alpha}\big)
\partial_2 \udel \eta^{2l}\dx \geq 0 \psp .
\]
The right-hand side of \reff{apri 7} is estimated with the help of the Cauchy-Schwarz inequality
which gives  for any $\eps >0$
\begin{eqnarray*}
\lefteqn{\iom D^2f_\delta(\nabla \udel) (\nabla \partial_2 \udel, \nabla \eta)\eta^{2l-1} 
\gam{2}{\alpha} \partial_2\udel \dx}\\
&\leq& \eps \iom D^2f_\delta(\nabla \udel) (\nabla \partial_2\udel ,\nabla \partial_2 \udel) 
\eta^{2l} \gam{2}{\alpha}\dx\\ 
&&+ c(\eps) \iom D^2f_\delta(\nabla \udel) (\nabla \eta,\nabla \eta) 
\eta^{2l-2}\gam{2}{\alpha} |\partial_2 \udel|^2 \dx
\end{eqnarray*}
and we have \reff{apri 4} by choosing $\eps$ sufficiently small.\\

We estimate (recall \reff{apri 1} and \reff{apri 2})
\begin{eqnarray}\label{apri 8}
\lefteqn{\iom D^2 f_\delta(\nabla \udel)\big(\nabla \eta,\nabla \eta\big) \eta^{2l-2} \gam{2}{\alpha+1}\dx}\nonumber\\
&\leq & c \Bigg[ \iom |\nabla \eta|^2 \eta^{2l-2}\gam{1}{\frac{\gamma}{2}}\gam{2}{\alpha+1}\dx 
+ \iom |\nabla \eta|^2 \eta^{2l-2} \gam{2}{\alpha+1+\frac{p-2}{2}}\dx\Bigg]\psp .
\end{eqnarray}
Let $\gamma >0$ and define the numbers $p_1$, $p_2$ via
\[
1 < p_1= \frac{1}{\gamma}\psp , \fsp p_2 = \frac{1}{1-\gamma} \psp .
\]
Using Young's inequality we obtain for the first integral on the right-hand side of \reff{apri 8}
\begin{eqnarray}\label{apri 9}
\iom |\nabla \eta|^2 \eta^{2l-2}\gam{1}{\frac{\gamma}{2}} \gam{2}{\alpha+1}\dx
&=& \iom \Big[|\nabla \eta|^2 \eta^{2l-2}\Big]^{\frac{1}{p_1}} \gam{1}{\frac{\gamma}{2}} 
\Big[|\nabla \eta|^2 \eta^{2l-2}\Big]^{\frac{1}{p_2}} \gam{2}{\alpha+1}\dx\nonumber\\
&\leq & c \Bigg[ 1+ \iom |\nabla \eta|^2 \eta^{2l-2} \gam{2}{\frac{\alpha+1}{1-\gamma}} \dx \Bigg] \psp .
\end{eqnarray}
With \reff{apri 8} and \reff{apri 9} the proof of Proposition \ref{apri prop 1} is finished by observing
\begin{eqnarray*}
\lefteqn{\iom \eta^{2l} \gam{2}{\alpha + \frac{p-2}{2}} |\partial_{22}\udel|^2 \dx}\nonumber\\
&\leq& c \iom  f_{2}''(\partial_2\udel)|\partial_{22}\udel|^2\eta^{2l} \gam{2}{\alpha} \dx \nonumber\\
&\leq& c \iom \Bigg[ f_{1,\delta}''(\partial_1\udel)|\partial_{12}\udel|^2 
+ f''_{2,\delta}(\partial_2 \udel) |\partial_{22} \udel|^2\Bigg]\eta^{2l} \gam{2}{\alpha}\dx\nonumber\\
&\leq & c \iom D^2f_\delta(\nabla \partial_2\udel, \nabla \partial_2 \udel) \eta^{2l} \gam{2}{\alpha}\dx \psp . \qed
\end{eqnarray*}

\emph{Step 3. Main inequality.}
\begin{proposition}\label{apri prop 2}
Given the above hypotheses let for some $\tau_s >0$, $\tau_\alpha >0$
\[
0 \leq s := \frac{p-2}{2}+\tau_s \psp , \fsp  \alpha := -\frac{1}{2} + \tau_\alpha \psp . 
\]
Then for $l$ sufficiently large and a local constant $c(\eta,l)$ independent of $\delta$, it holds
\begin{eqnarray}\label{apri 10}
\iom \eta^{2l} \gam{2}{s+1} \dx &\leq& c \Bigg[1
 + \iom |\nabla \eta|^2 \eta^{2l-2}\gam{2}{\frac{\alpha+1}{1-\gamma}} \dx
+ \iom |\nabla \eta|^2 \eta^{2l-2}\gam{2}{\alpha +1+\frac{p-2}{2}} \dx\nonumber\\ 
&&+ \iom \gam{2}{2s - \alpha -  \frac{p-2}{2}}\eta^{2l} \dx \Bigg] \psp .
\end{eqnarray}
\end{proposition}

\emph{Proof of Proposition \ref{apri prop 2}.} Let us first note that on the left-hand side of \reff{apri 10}
we have $s+1 = (p/2) + \tau_s$. Moreover, since $\|\udel\|_{L^\infty(\Omega)} \leq c$, we estimate
\begin{eqnarray*}
\iom |\partial_2 \udel|^2 \gam{2}{s}\eta^{2l} \dx &= &
\iom \partial_2 \udel \partial_2 \udel \gam{2}{s} \eta^{2l} \dx
= - \iom \udel \partial_2 \Big[\partial_2 \udel \gam{2}{s}\eta^{2l}\Big]\dx\\
&\leq & c \Bigg[ \iom |\partial_{22} \udel| \gam{2}{s} \eta^{2l} \dx +
\iom |\partial_2 \udel| \eta^{2l-1} |\nabla \eta|\gam{2}{s}\dx\\
&&+  \iom \gam{2}{s-1}|\partial_2 \udel|^2 |\partial_{22} \udel| \eta^{2l} \dx \Bigg]\\
&\leq& c \Bigg[ \iom |\partial_{22}\udel| \gam{2}{s} \eta^{2l} \dx \\
&&+ \eps \iom |\partial_2 \udel|^2 \gam{2}{s} \eta^{2l} \dx 
+ c(\eps) \iom |\nabla \eta|^2 \eta^{2l-2}\gam{2}{s} \dx \Bigg]\psp ,
\end{eqnarray*}
leading to ($\eps >0$ sufficiently small)
\begin{equation}\label{apri 11}
\iom |\partial_2 \udel|^2 \gam{2}{s} \eta^{2l} \dx \leq
c \Bigg[\iom |\partial_{22} \udel| \gam{2}{s} \eta^{2l} \dx +
\iom |\nabla \eta|^2 \eta^{2l-2}\gam{2}{s}\dx \Bigg] \psp .
\end{equation}

The first term on the right-hand side of \reff{apri 11} is estimated with the help of Young's inequality
\begin{eqnarray}\label{apri 12}
\lefteqn{\iom |\partial_{22} \udel| \gam{2}{s} \eta^{2l} \dx =
\iom |\partial_{22} \udel|\gam{2}{\frac{\alpha}{2}+\frac{p-2}{4}}
\gam{2}{s-\frac{\alpha}{2} - \frac{p-2}{4}} \eta^{2l} \dx}\nonumber\\
&\leq & c\Bigg[\iom |\partial_{22} \udel|^2 \gam{2}{\alpha + \frac{p-2}{2}} \eta^{2l} \dx
+\iom \gam{2}{2s-\alpha - \frac{p-2}{2}} \eta^{2l} \dx\Bigg] \psp .
\end{eqnarray}
Now \reff{apri 5} is applied to the first term on the right-hand side of \reff{apri 12} which gives
\begin{eqnarray}\label{apri 13}
\lefteqn{\iom |\partial_{22} \udel| \gam{2}{s} \eta^{2l} \dx}\nonumber\\ 
&\leq&
c \Bigg[ 1 + \iom |\nabla \eta|^2 \eta^{2l-2}\gam{2}{\frac{\alpha+1}{1-\gamma}}\dx
+ \iom |\nabla \eta|^2 \eta^{2l-2}\gam{2}{\alpha +1+\frac{p-2}{2}} \dx\nonumber\\ 
&&+ \iom \gam{2}{2s - \alpha -  \frac{p-2}{2}}\eta^{2l} \dx \Bigg] \psp .
\end{eqnarray}

Combining \reff{apri 11} and \reff{apri 13} we finally obtain
\begin{eqnarray}\label{apri 14}
\iom \gam{2}{s+1} \eta^{2l}\dx &\leq & c \Bigg[ 1 + \iom \big(\eta^{2l} + |\nabla \eta|^2 \eta^{2l-2}\big) \gam{2}{s} \dx\nonumber\\
&& + \iom |\nabla \eta|^2 \eta^{2l-2}\gam{2}{\frac{\alpha+1}{1-\gamma}}
+ \iom |\nabla \eta|^2 \eta^{2l-2}\gam{2}{\alpha +1+\frac{p-2}{2}} \dx\nonumber\\ 
&&+ \iom \gam{2}{2s - \alpha -  \frac{p-2}{2}}\eta^{2l} \dx \Bigg] \psp .
\end{eqnarray}
The first intergral on the right-hand side of \reff{apri 14} can be absorbed in the left-hand side whenever $l$ is sufficiently
large (compare, e.g., \cite{BF:2020_3}, Proof of Proposition 2.2) which completes the proof of the proposition. \qed\\

\emph{Step 4. Conclusion.}
To finish the proof of Theorem \ref{apri theo 1} we observe that we may exactly follow the lines of
\cite{BF:2020_3}, Theorem 1.3,  provided that
\begin{equation}\label{apri 15}
\alpha + 1 + \frac{p-2}{2} = \alpha + \frac{p}{2} < s+1 \psp , \fsp
2s - \alpha - \frac{p-2}{2} < s + 1\psp ,
\end{equation}
and provided that
\begin{equation}\label{apri 16}
\alpha + 1 < (s+1)(1-\gamma) = s+1 - (s+1) \gamma \psp .
\end{equation}
We note that \reff{apri 15} is equivalent to
\begin{equation}\label{apri 17}
\tau_\alpha < \tau_s + \frac{1}{2}\psp , \fsp \tau_s < \tau_\alpha + \frac{1}{2}\psp ,
\fsp\mbox{i.e.}\msp |\tau_s -\tau_\alpha| < \frac{1}{2}\psp .
\end{equation}
Condition \reff{apri 16} turns into
\[
-\frac{1}{2}+\tau_\alpha < \frac{p-2}{2} + \tau_s - \frac{p+2\tau_s}{2} \gamma \psp ,
\]
which can be written as
\begin{equation}\label{apri 18}
\gamma < \frac{p-1 + 2 (\tau_s-\tau_\alpha)}{p+2\tau_s} \psp .
\end{equation}

\begin{enumerate}
\item If $\gamma = 0$, then \reff{apri 17} and \reff{apri 18} are satisfied for any $\tau_s > 1/4$ and with the choice 
$\tau_\alpha = \tau_s - 1/4$, which implies the first claim of the theorem.
\item If we have \reff{apri 3} with $\gamma >0$, then we choose $\tau_\alpha > 0$ sufficiently small and 
$\tau_s = \frac{1}{2} + \frac{\tau_\alpha}{2}$. Then we have \reff{apri 17} and \reff{apri 18} for $\tau_\alpha$
sufficiently small on account of our assumption \reff{apri 3}. 
We note that with the notation introduced above we have
\[
\chi = 2 (s+1) = 2 \Bigg[\frac{p-2}{2} + \tau_s +1 \Bigg] = p + 2\tau_s \psp .
\] 
The choice $\tau_s > 1/2$ completes the proof of Theorem \ref{apri theo 1}. \hspace*{\fill}\qed
\end{enumerate} 

In our final theorem on regularity we pose the question, whether am improved ellipticity for the linear part is compatible
with $p$-growth for arbitrary $p$ in the following sense: do we have higher integrability of $\partial_2 u$
and simulaneously of $\partial_1 u$? In fact, in addition to Theorem \ref{apri theo 1} we have

\begin{theorem}\label{apri theo 2}
Suppose that we have \reff{intro not 1} and the first condition of \reff{intro not 2}, i.e.~a linear growth condition for
$f_1$. Moreover, suppose that we have \reff{intro ell 1} with $2 < \mu$ and $\gamma =0$ and let assumption
\reff{apri 2} hold.\\

Then there exists $u\in \mathcal{M}$ such that for any $\kappa < \infty$
\[
\partial_1 u \in L^\kappa(\Omega) \psp .
\]
\end{theorem}

\emph{Proof of Theorem \ref{apri theo 2}.} 
Since we already discussed Theorem \ref{apri 1818 theo 1} covering the case $p=2$  we may suppose that $p >2$.
Given the regularization of Theorem \ref{apri theo 1} we need a counterpart for
Propostion \ref{apri prop 1} with $\partial_{22} \udel$ rplaced by $\partial_{11}\udel$ and $\gam{2}{}$ by $\gam{1}{}$.

\begin{proposition}\label{apri prop 3}
Given the hypotheses of Theorem \ref{apri theo 2} we have for any real numbers $\chi > 1$, $\alpha > -1/2$
\begin{eqnarray}\label{apri lin 1}
\lefteqn{\iom \eta^{2l} \gam{1}{\alpha - \frac{\mu}{2}} |\partial_{11}\udel|^2 \dx}\nonumber\\
&\leq& c \Bigg[ 1+ \iom |\nabla \eta|^2 \eta^{2l-2}\gam{1}{\alpha+1}\dx
+ \iom |\nabla \eta|^2 \eta^{2l-2}\gam{1}{(\alpha +1)\frac{\chi}{\chi-(p-2)}} \dx \Bigg] \psp .
\end{eqnarray}
\end{proposition}

\emph{Proof of Proposition \ref{apri prop 3}}. The counterpart of \reff{apri 8} reads as
\begin{eqnarray}\label{apri lin 2}
\lefteqn{\iom D^2 f_\delta(\nabla \udel)\big(\nabla \eta,\nabla \eta\big) \eta^{2l-2} \gam{1}{\alpha+1}\dx}\nonumber\\
&\leq & c \Bigg[ \iom |\nabla \eta|^2 \eta^{2l-2} \gam{1}{\alpha+1}\dx 
+ \iom |\nabla \eta|^2 \eta^{2l-2} \gam{1}{\alpha+1}\gam{2}{\frac{p-2}{2}}\dx\Bigg]\psp .
\end{eqnarray}
Now we let (recall $p >2$)
\[
p_1 = \frac{\chi}{\chi - (p-2)} \psp , \fsp
1 < p_2= \frac{\chi}{p-2}\psp .
\]
For the choice of $\chi$ we have Theorem \ref{apri theo 1}, $i$), in mind.
Now Young's inequality implies for the second integral on the right-hand side of \reff{apri lin 2}
\begin{eqnarray}\label{apri lin 3}
\iom |\nabla \eta|^2 \eta^{2l-2}\gam{1}{\alpha+1}\gam{2}{\frac{p-2}{2}}\dx
&=& \iom \Big[|\nabla \eta|^2 \eta^{2l-2}\Big]^{\frac{1}{p_1}} \gam{1}{\alpha +1} 
\Big[|\nabla \eta|^2 \eta^{2l-2}\Big]^{\frac{1}{p_2}} \gam{2}{\frac{p-2}{2}}\dx\nonumber\\
&\leq & c \Bigg[ 1+ \iom |\nabla \eta|^2 \eta^{2l-2} \gam{1}{(\alpha+1) \frac{\chi}{\chi-(p-2)}} \dx \Bigg] \psp .
\end{eqnarray}
By \reff{apri lin 2} and \reff{apri lin 3} we have \reff{apri lin 1}, i.e.~the proposition is proved. \qed\\

Now a variant of Proposition  \ref{apri prop 1} is needed.

\begin{proposition}\label{apri prop 4}
Given the hypotheses of Theorem \ref{apri theo 2} let for some $\tau_s >0$, $\tau_\alpha >0$
\[
0 \leq s := -\frac{1}{2}+\tau_s \psp , \fsp  \alpha := -\frac{1}{2} + \tau_\alpha \psp . 
\]
Then for $l$ sufficiently large and a local constant $c(\eta,l)$ independent of $\delta$ it holds
\begin{eqnarray}\label{apri lin 4}
\iom \eta^{2l} \gam{1}{s+1} \dx &\leq& c \Bigg[1
 + \iom |\nabla \eta|^2 \eta^{2l-2}\gam{1}{\alpha+1}\dx\nonumber\\
&&+ \iom |\nabla \eta|^2 \eta^{2l-2}\gam{1}{(\alpha +1)\frac{\chi}{\chi-(p-2)}} \dx\nonumber\\
&&+ \iom \gam{1}{2s - \alpha + \frac{\mu}{2}}\eta^{2l} \dx \Bigg] \psp .
\end{eqnarray}
\end{proposition}

\emph{Proof of Proposition \ref{apri prop 4}.} Here we have instead of \reff{apri 11}, \reff{apri 12}
\begin{equation}\label{apri lin 5}
\iom |\partial_1 \udel|^2 \gam{1}{s} \eta^{2l} \dx \leq
c \Bigg[\iom |\partial_{11} \udel| \gam{1}{s} \eta^{2l} \dx +
\iom |\nabla \eta|^2 \eta^{2l-2}\gam{1}{s}\dx \Bigg] 
\end{equation}
and
\begin{eqnarray}\label{apri lin 6}
\lefteqn{\iom |\partial_{11} \udel| \gam{1}{s} \eta^{2l} \dx =
\iom |\partial_{11} \udel|\gam{1}{\frac{\alpha}{2}-\frac{\mu}{4}}
\gam{1}{s-\frac{\alpha}{2} + \frac{\mu}{4}} \eta^{2l} \dx}\nonumber\\
&\leq & c\Bigg[\iom |\partial_{11} \udel|^2 \gam{1}{\alpha - \frac{\mu}{2}} \eta^{2l} \dx
+\iom \gam{1}{2s-\alpha + \frac{\mu}{2}} \eta^{2l} \dx\Bigg] \psp .
\end{eqnarray}
We then obtain
\begin{eqnarray}\label{apri lin 7}
\iom \gam{1}{s+1} \eta^{2l}\dx &\leq & c \Bigg[ 1 + \iom \big(\eta^{2l} + |\nabla \eta|^2 \eta^{2l-2}\big) \gam{1}{s} \dx\nonumber\\
&& + \iom |\nabla \eta|^2 \eta^{2l-2}\gam{1}{\alpha+1}
+ \iom |\nabla \eta|^2 \eta^{2l-2}\gam{1}{(\alpha +1)\frac{\chi}{\chi-(p-2)}} \dx\nonumber\\ 
&&+ \iom \gam{1}{2s - \alpha + \frac{\mu}{2}}\eta^{2l} \dx \Bigg] \psp ,
\end{eqnarray}
and \reff{apri lin 7} gives the proposition. \hspace*{\fill}\qed \\

Finally we arrive at the following choice of parameters. 
\begin{eqnarray}\label{apri lin 8}
\alpha + 1 < s+1 & \Leftrightarrow & \tau_\alpha < \tau_ s \psp ,\\[2ex]
2s - \alpha + \frac{\mu}{2} < s+1 &\Leftrightarrow& \tau_s - \tau_\alpha < 1 - \frac{\mu}{2}\psp .\label{apri lin 9}
\end{eqnarray}
Moreover, the inequality
\begin{equation}\label{apri lin 10}
\Bigg[\frac{1}{2}+\tau_\alpha\Bigg] \psp \Bigg[\frac{\chi}{\chi-(p-2)}\Bigg] < \frac{1}{2} + \tau_s
\end{equation}
needs to be true.\\

For any given $\tau_s >0$ we may choose $\tau_\alpha >0$ such that \reff{apri lin 8} and \reff{apri lin 9} hold (recalling
$\mu < 2$). On account of Theorem \ref{apri theo 1} we then choose $\chi$ sufficiently large such that
\reff{apri lin 10} is satisfied as well. This proves the theorem. \hspace*{\fill}\qed \\

\begin{remark}\label{apri rem 2}
We note that even higher integrabilty of $\partial_2 u$ with some $\chi > p+1$ (compare Theorem \ref{apri theo 1}, $ii$)) 
would give some higher integrability of $\partial_1 u$ by choosing 
\[
\tau_\alpha := \frac{3(2-\mu) - (p-2)}{2 (p-2)}
\]
whenever $3(2-\mu)-(p-2) > 0$. This indicates that we also may prove variants of Theorem \ref{apri theo 2} with some
weaker hypotheses. This is left to the reader.
\end{remark}

Finally we observe that with the higher integrability of the gradient w.r.t.~any exponent we have using the idea, e.g., of \cite{Bi:1818},
Section 4.2.4 we obtain as a corollary:

\begin{cor}\label{apri cor 1}
Suppose that we have the hypotheses of Theorem \ref{apri theo 2}. Then any weak $L^1$-cluster-point of the regularizing sequence
$\{\udel\}$ is of class $C^{1,\alpha}(\Omega)$.
\end{cor}
\section{Uniqueness of solutions}\label{uni}

The first uniqueness result is already established in Theorem \ref{apri 1818 theo 1} and Corollary
\ref{apri 1818 cor 1}.\\

We now finally discuss the case that the superlinear part of the energy density under consideration
is given in terms of a $N$-function.

In this case the uniqueness of solutions to the generalized problems can be established without using
the results of the previous section.\\

The reason is the existence of an unique dual solution which, by the duality relation, can be carried over to 
establish the uniqueness of generalized minimizers.\\

So let us suppose that we have the situation as described in the hypotheses of Theorem \ref{dual theo 1}, in particular
\[
\sigma: = Df\big(\nabla^a \ovu\big) = \Big(f_1'\big(\partial_1^a \ovu\big), A'\big(|\partial_2 \ovu|\big)\Big)
\]
is a solution of the dual variational problem \reff{dual 10} whenever $\ovu$ is a generalized minimizer in the sense
of  \reff{probbv}.\\

Now, if we have a closer look at the arguments from measure theory leading to Theorem 7 of \cite{Bi:2000_1},
then we may adapt the proof to the situation at hand and obtain:

\begin{theorem}\label{uni dual theo 1}
Suppose that we have the hypotheses stated in the beginning of Section \ref{apriori p} in front of Theorem \ref{apri theo 1}. 
Moreover, suppose
that $u^*\in \mathcal{M}$ is a $L^1$-cluster-point of a sequence $\{\udel\}$ from $u_0+\wnull$ such that for some $\varkappa \in \rz$
the functions
\begin{equation}\label{dual uni 0.0}
h_\delta := h^1_\delta+h^2_\delta := \big(1+|\partial_ 1\udel|^2\big)^{\frac{\varkappa}{2}} +
\big(1+|\partial_ 2\udel|^2\big)^{\frac{\varkappa}{2}}
\end{equation}
are of class $W^1_{1,\op{loc}}(\Omega)$ uniform w.r.t.~$\delta \in (0,1)$.\\

Then the dual problem
\reff{dual 10} is uniquely solvable. 
\end{theorem}

We emphasize that Theorem \ref{uni dual theo 1} is valid without any restriction on the exponent $\mu$ in \reff{intro ell 1}.
This is just needed for Corollary \ref{uni dual cor 1} below.\\

\emph{Main idea of the proof of Theorem \ref{uni dual theo 1}.} 
Using the notation of \cite{Bi:1818}, Section 2.2, we first note that the unboundedness of $U:=\op{Im}(\nabla f)$ w.r.t.~the direction $\xi_2$
causes no essential changes in the previous arguments.\\

Given the regularization $\{\udel\}$ of Step 1 in the proof of Theorem \ref{apri theo 1}, the main changes concern the proof of
(compare p.~22 of \cite{Bi:1818}) 
\begin{equation}\label{dual uni 0.1}
\sigma_\delta(x) \to \sigma(x)\psp ,\fsp \delta |\partial_1 \udel(x)|^{p-1} \to 0 \psp , 
\end{equation}
on $\Sigma$ as $\delta \to 0$. In \cite{Bi:1818}, these convergences follow form the uniform local $W^1_2$-regularity of 
$\sigma_\delta$ which in \cite{Bi:1818} is a consequence of
\begin{equation}\label{dual uni 0.2}
D^2f(\xi)\big(\eta,\eta) \leq \big(1+|\xi|\big)^{-1} |\eta|^2\psp ,\fsp \xi ,\psp \eta \in \rz^2 \psp .
\end{equation}
Since \reff{dual uni 0.2} in general is no longer valid, we use \reff{dual uni 0.0} to get
the a.e.~convergence of $\nabla \udel$. \qed\\

We finally have
\begin{cor}\label{uni dual cor 1}
Suppose that we have the assumptions of Theorem \ref{apri theo 2} with $p >2$.
Then we have the uniqueness of generalized solutions of problem \reff{probbv}.
\end{cor}

\emph{Proof of Corollary \ref{uni dual cor 1}.} We first establish the uniqueness of the dual solution $\sigma$.
To this purpose we recall \reff{apri 5} and \reff{apri lin 1}, i.e.~given the regularization introduced in Step 1 of the proof of 
Theorem \ref{apri theo 1} we have choosing $\alpha =0$
\begin{eqnarray}\label{uni dual 1}
\lefteqn{\iom D^2 f_\delta (\nabla \udel) \big(\nabla \partial_1 \udel,\nabla \partial_1\udel\big) \eta^{2} \dx}\nonumber\\
&\leq&  c \Bigg[ 1+ \iom |\nabla \eta|^2 \gam{1}{}\dx
+ \iom |\nabla \eta|^2 \gam{1}{\frac{\chi}{\chi-(p-2)}} \dx \Bigg] 
\end{eqnarray}
as well as
\begin{eqnarray}\label{uni dual 2}
\lefteqn{\iom D^2 f_\delta (\nabla \udel) \big(\nabla \partial_2 \udel,\nabla \partial_2\udel\big) \eta^{2}\dx}\nonumber\\
&\leq &  c \Bigg[ 1+ \iom |\nabla \eta|^2 \gam{2}{}\dx
+ \iom |\nabla \eta|^2 \gam{2}{\frac{p}{2}} \dx \Bigg] \psp .
\end{eqnarray}
Now we compute
\begin{eqnarray*}
\iom \eta^{2}|\nabla h_\delta| \dx &\leq & c \iom \eta^{2l} \big(1+|\partial_1 u|^2\big)^{\frac{\varkappa-1}{2}} |\nabla \partial_1 \udel|\dx\\
 && +c \iom \eta^{2} \big(1+|\partial_2 u|^2\big)^{\frac{\varkappa -1}{2}} |\nabla \partial_2 \udel|\dx\\
&\leq &c \iom\eta^{2} \Big[\gam{1}{-\frac{\mu}{4}} |\partial_1 \partial_1 \udel| + |\partial_2 \partial_1 \udel|\Big]
\gam{1}{\frac{\mu}{4}+\frac{\varkappa-1}{2}}\dx\\
& & +c \iom\eta^{2} \Big[\gam{1}{-\frac{\mu}{4}} |\partial_1 \partial_2 \udel| + |\partial_2 \partial_2 \udel|\Big]
\gam{1}{\frac{\mu}{4}}\gam{2}{\frac{\varkappa-1}{2}}\dx
\end{eqnarray*}
which leads to 
\begin{eqnarray}\label{uni dual 3}
\iom \eta^{2}|\nabla h_\delta| \dx 
&\leq & \iom \eta^{2}D^2f_\delta(\nabla \udel)\big(\nabla \partial_1 \udel, \nabla \partial_1 \udel\big)\dx\nonumber\\
&&+ \iom \eta^{2}D^2f_\delta(\nabla \udel)\big(\nabla \partial_2 \udel, \nabla \partial_2 \udel\big)\dx\nonumber\\
&&+ \iom \eta^2  \gam{1}{t_1}\dx + \iom \eta^2 \gam{2}{t_2}\dx
\end{eqnarray}
with some finite exponents $t_1$ and $t_2$. By \reff{uni dual 1} - \reff{uni dual 3} and Theorem \ref{apri theo 1}, $i$), Theorem \ref{apri theo 2}
we have that $h_\delta$ is uniformly of class $W^{1}_{1,\op{loc}}$ which by Theorem \ref{uni dual theo 1} implies the
uniqueness of the dual solution.\\

The next step for proving Corollary \ref{uni dual cor 1} is to use Corollary \ref{apri cor 1} to obtain $C^{0,\alpha}$-regularity
of the dual solution. Then, as outlined in the proof of Theorem A.9 in \cite{Bi:1818}, a suitable comparison argument w.r.t.~$\sigma$
can be carried out to obtain for any generalized minimizer $\ovu \in \cbv$
\begin{equation}
\nabla \ovu = \nabla f^*(\sigma) \psp .
\end{equation}
We note that \reff{uni dual 1} in particular gives
\[
\nabla^s \ovu = 0 \psp .
\]
By the uniqueness of $\sigma$ the uniqueness of generalized minimizers up to an additive constant is established. Finally,
on account ot $\ovu  \in \cbv$ we have $(\ovu -u_0)\nu_2= 0$ $\mathcal{H}^1$ a.e.~on $\partial \Omega$ and in conclusion
the uniqueness of generalized minimizers. \qed\\



\begin{tabular}{ll}
Michael Bildhauer&bibi@math.uni-sb.de\\
Martin Fuchs&fuchs@math.uni-sb.de\\[5ex]
Department of Mathematics&\\
Saarland University&\\
P.O.~Box 15 11 50&\\
66041 Saarbr\"ucken&\\ 
Germany&
\end{tabular}

\end{document}